\documentclass[a4paper, 11pt, final]{article}

\usepackage[all]{xy}
\usepackage{style}
 \usepackage{amsthm, eufrak}
\usepackage{dsfont}
 \usepackage{geometry}         

\title{Operations in Milnor $K$-theory} 
                                
\author{Charles Vial} 
\date{} 

\begin{document}

\maketitle                 

\begin{abstract} We show that operations in Milnor $K$-theory mod $p$ of
  a field are spanned by divided power operations. After giving an
  explicit formula for divided power operations and extending them to
  some new cases, we determine for all fields $k_0$ and all prime
  numbers $p$, all the operations $\KM_i/p \r \KM_j/p$ commuting with
  field extensions over the base field $k_0$. Moreover, the integral
  case is discussed and we determine the operations $\KM_i/p \r
  \KM_j/p$ for smooth schemes over a field.

\end{abstract}


\section*{Introduction}

Let $k_0$ be a field and $p$ be a prime number different from the
characteristic of $k_0$. In \cite{Voevodsky1}, Voevodsky constructs
Steenrod operations on the motivic cohomology $H^{*,*}(X,\Z/p)$ of a
general scheme over $k_0$. However, when $p$ is odd or when $p=2$ and
$-1$ is a square in $k_0^\times$, such operations vanish on the
motivic cohomology groups $H^{i,i}(\mathrm{Spec} \ k,\Z/p)$ for $i>0$
of the spectrum of a field extension $k$ of $k_0$. Here, we study
operations on $H^{i,i}(\mathrm{Spec} \ k,\Z/p)$ which are defined only
for fields.

The same phenomenon happens in \'etale cohomology, where Steenrod
operations, as defined by Epstein in \cite{Epstein}, vanish on the
\'etale cohomology $H^i_{et}(\mathrm{Spec} \ k,\Z/p)$ of a field if
$p$ is odd or if $p=2$ and $\sqrt{-1} \in k$.  Under the assumption of
the Bloch-Kato conjecture, our operations give secondary operations
relatively to Steenrod operations on the \'etale cohomology of fields.
\medskip

Given a base field $k_0$ and a prime number $p$, an operation on
$\KM_i /p$ is a function $\KM_i(k)/p \r \KM_*(k)/p$ defined for all
fields $k / k_0$, compatible with extension of fields. In other words,
it is a natural transformation from the functor $ \KM_i/p :
\mathsf{Fields}_{/k_0} \r \mathsf{Sets}$ to the functor $ \KM_*/p :
\mathsf{Fields}_{/k_0} \r \mathbf{F}_p - \mathsf{Algebras}$. It is
important for our purpose that our operations should be functions and
not only additive functions, the reason being that additive operations
will appear to be trivial in some sense (see section 3.5). In these
notes, we determine all operations $\KM_i/p \r \KM_*/p$ over any field
$k_0$, no matter if $p \neq \mathrm{char} \ k_0$ or not. This is
striking, especially in the case when $i=1$.\medskip

Let $n$ be a non-negative integer and $k$ any field. Let $x =
\sum_{r=1}^l s_r$ be a sum of $l$ symbols in $\KM_i(k)/p$, the mod $p$
Milnor $K$-group of $k$ of degree $i$. We define the $n^\mathrm{th}$
divided power of $x$, given as a sum of symbols, by $$\gamma_n(x)=
\sum_{1 \leq l_1 < \ldots < l_n \leq l} s_{l_1} \cdot \ldots \cdot
s_{l_n} \ \in \KM_{ni}(k)/p.$$ Such a divided power may depend on the
way $x$ has been written as a sum of symbols and thus a well-defined
map $\gamma_n : \KM_i(k)/p \rightarrow \KM_{ni}(k)/p$ may not exist.
However, $\gamma_0(x) =1$ and $\gamma_1(x) = x$ and as such,
$\gamma_0$ and $\gamma_1$ are always well-defined. The axioms for
divided powers (see prop. \ref{prop}) formalize the properties of
$\frac{x^n}{n!}$ in a $\Q$-algebra, see \cite{BO} for some general
discussion of a divided power structure on an ideal in a commutative
ring. In his paper \cite{Kafields}, Kahn shows that the above formula
gives well-defined divided powers $\gamma_n : \KM_{2i}(k)/p \r
\KM_{2ni}(k)/p$ for $p$ odd and $\gamma_n : \KM_i(k)/2 \r
\KM_{ni}(k)/2$ for $k$ containing a square root of $-1$. Kahn's result
is based on previous work by Revoy on divided power algebras,
\cite{Revoy}. Divided powers are also mentioned in a letter of Rost to
Serre, \cite{Galois}. In this paper, we show that in these cases,
divided powers define operations in the above sense and form a basis
for all possible operations on mod $p$ Milnor $K$-theory.

In the remaining case, when $-1$ is not a square in the base field
$k_0$, divided powers as defined above are not well-defined on mod $2$
Milnor $K$-theory. However, we will define some new, weaker
operations, and show that these new operations are all the possible
operations on mod $2$ Milnor $K$-theory. \medskip

Precisely, we will prove :

\begin{theorem2}[$p$ odd] \label{main} Let $k_0$ be any field, $p$ an
  odd prime number. The algebra of operations $\KM_i(k)/p \r
  \KM_*(k)/p$ commuting with field extensions over $k_0$ is
  \begin{itemize}
  \item If $i=0$, the free $\KM_*(k_0)/p$-module of rank $p$ of
    functions $\mathbf{F}_p \r \KM_*(k_0)/p$.
  \item If $i \geq 1$ is odd, the free $\KM_*(k_0)/p$-module
    $$\KM_*(k_0)/p\cdot \gamma_0 \ \oplus \ \KM_*(k_0)/p\cdot
    \gamma_1.$$
  \item If $i \geq 2$ is even, the free $\KM_*(k_0)/p$-module $$
    \bigoplus_{n\geq 0} \KM_*(k_0)/p \cdot \gamma_n.$$

  \end{itemize}
\end{theorem2}

\begin{theorem2} [$p=2$] \label{main2} Let $k_0$ be any field. The
  algebra of operations $\KM_i(k)/2 \r \KM_*(k)/2$ commuting with
  field extensions over $k_0$ is
\begin{itemize}
\item If $i=0$, the free $\KM_*(k_0)/2$-module of rank $2$ of
  functions $\mathbf{F}_2 \r \KM_*(k_0)/2$.
\item If $i=1$, the free $\KM_*(k_0)/2$-module of rank $2$, generated
  by $\gamma_0$ and $\gamma_1$.
\item If $i \geq 2$, the $\KM_*(k_0)/2$-module $$\KM_*(k_0)/2\cdot
  \gamma_0 \oplus \KM_*(k_0)/2\cdot \gamma_1 \oplus \bigoplus_{n\geq
    2} \mathrm{Ker}(\tau_i) \cdot \gamma_n,$$ where $\tau_i :
  \KM_*(k_0)/2 \r \KM_*(k_0)/2$ is the map $\ x \mapsto \{-1\}^{i-1}
  \cdot x$.

\end{itemize}
 \end{theorem2}

 Actually, the divided powers $\gamma_n$ are not always defined with
 the assumptions of Theorem \ref{main2}.  However, if $y_n \in
 \mathrm{Ker}(\tau_i)$, the map $y_n \cdot \gamma_n$ will be shown to
 be well-defined. Notice that when $-1$ is a square in $k_0$, the map
 $\tau_i$ is the zero map, and hence $\mathrm{Ker}(\tau_i) =
 \KM_*(k_0)/2$.

 Also, the divided powers satisfy the relation $\gamma_m \cdot
 \gamma_n = \begin{pmatrix} m+n \\ n
 \end{pmatrix} \gamma_{m+n}$. Together with the algebra structure on
 $\KM_*(k_0)/p$, this gives the algebra structure of the algebra of
 operations $\KM_i/p \r \KM_*/p$ over $k_0$. In fact, divided powers
 satisfy all the relations mentionned in properties \ref{prop} and
 make Milnor $K$-theory into a divided power algebra in Revoy's
 notation \cite{Revoy}. \medskip

 As Nesterenko-Suslin \cite{NS} and Totaro \cite{Totaro} have shown,
 there is an isomorphism \\ $H^{n,n}(\mathrm{Spec} \ k, \Z)
 \stackrel{\simeq}{\longrightarrow} \KM_n(k)$ where
 $H^{n,n}(\mathrm{Spec} \ k, \Z)$ denotes motivic cohomology. This
 isomorphism, together with Theorem \ref{main}, provides operations on
 the motivic cohomology groups $H^{n,n}(\mathrm{Spec} \ k, \Z/p)$.
 Also, since the Bloch-Kato conjecture seems to have been proven by
 Rost and Voevodsky (see \cite{Voevodsky2} and Weibel's paper
 \cite{Weibel} that patches the overall proof by using operations from
 integral cohomology to $\Z/p$ cohomology avoiding lemma 2.2 of
 \cite{Voevodsky2} which seems to be false as stated), this gives the
 operations in Galois cohomology of fields, with suitable
 coefficients.  \medskip

 We also describe some new operations in integral Milnor $K$-theory over
 any base field $k_0$. Under some reasonable hypothesis on an
 operation $\varphi : \KM_i \r \KM_*$ defined over $k_0$, we are able
 to show that $\varphi$ is in the $\KM_*(k_0)$-span of our weak
 divided power operations. See sections 2.5 and 3.4. \medskip

 Finally, we are able to determine operations $\KM_i/p \r \KM_j/p$ in
 the more general setup of smooth schemes over a field $k$. The
 Milnor $K$-theory ring of a smooth scheme $X$ over $k$ is defined to
 be the subring of the Milnor $K$-theory of the function field $k(X)$
 whose elements are unramified along all divisors on $X$, i.e. which
 vanish under all residue maps corresponding to codimension-$1$ points
 in $X$. An operation $\KM_i/p \r \KM_j/p$ over the smooth $k$-scheme
 $X$ is a function that is functorial with respect to morphisms of
 $X$-schemes (see section 4).  Once again, if $p$ is odd and $i \geq
 2$ is even, or if $p=2$ and $k$ contains a square-root of $-1$, we
 have

 \begin{theorem2}\label{scheme}
   Operations $\KM_i/p \r \KM_*/p$ over the smooth $k$-scheme $X$ are
   spanned as a $\KM_*(X)/p$-module by the divided power operations.
 \end{theorem2} 

Assuming the Bloch-Kato conjecture, we obtain in this way all the
operations for the unramified cohomology of smooth schemes over the
field $k$. \medskip

The paper is organised as follows. In the first section, we start by
recalling some general facts about Milnor $K$-theory, particularly the
existence of residue and specialization maps. In the second section,
we give a detailed account on divided power operations and extend the
results mentioned in \cite{Kafields} to the case $p=2$, $\sqrt{-1}
\notin k^\times$.  We also describe some weak divided power operations
for integral Milnor $K$-theory.  Some applications to cohomological
invariants are discussed. Section 3 contains the main results, Theorem
\ref{main} and Theorem \ref{main2} are proven (see Propositions
\ref{main1}, \ref{odd} and \ref{mod2}), the integral case is discussed
(Proposition \ref{integral}) and additive operations are determined
integrally (Proposition \ref{additive}). In section 4, we extend our
results to the case of the Milnor $K$-theory of smooth schemes over a
field $k$ and prove Theorem \ref{scheme}.

Finally, in the last section we mention that all the previous results
hold in the more general setup of operations from Milnor $K$-theory to
cycle modules with a ring structure as defined by Rost in \cite{Rost}.
In particular, this determines all the operations mod $p$ from Milnor
$K$-theory to Quillen $K$-theory of a field. \medskip

\textit{Acknowledgement.} I am deeply indebted to Burt Totaro for this
work. He suggested the problem to me and gave me some helpful advice
through his detailed remarks. This work is part of a Ph.D. thesis; the
ENS, the EPSRC and Trinity College (Cambridge) are acknowledged for
financial support.

\section{General Facts}

All the results in this section can be found in \cite{Gille}, Chapter
7.\medskip

Let $k$ be a field. The $n^\mathrm{th}$ Milnor $K$-group $\KM_n(k)$ is
the quotient of the $n$-fold tensor power $(k^\times)^{\otimes n}$ of
the multiplicative group $k^\times$ of the field $k$ by the relations
$a_1 \otimes \ldots \otimes a_n=0$ as soon as $a_i+a_j=1$ for some $1
\leq i < j \leq n$. We write $\{a_1, \ldots, a_n\}$ for the image of
$a_1 \otimes \ldots \otimes a_n$ in $\KM_n(k)$, such elements are
called symbols. Thus, elements in $\KM_n(k)$ are sums of symbols. The
relation $\{x,1-x\}=0$ in $\KM_2(k)$ is often referred to as the
Steinberg relation. In particular $\KM_0(k)=\Z$ and
$\KM_1(k)=k^\times$. This construction is functorial with respect to
field extension. There is a cup-product operation $\KM_n(k) \times
\KM_m(k) \r \KM_{n+m}(k)$ induced by the tensor product pairing
$(k^\times)^{\otimes n} \times (k^\times)^{\otimes m} \r
(k^\times)^{\otimes (m+n)}$. We write $\KM_*(k)$ for the direct sum
$\bigoplus_{n \geq 0} \KM_n(k)$. As a general fact, for any elements
$x$ and $y$ in $k^\times$, we have the relation $\{x,y\}=-\{y,x\}$.
Thus, cup product turns $\KM_*(k)$ into a graded commutative algebra.
We now state the easy but important

\begin{remark}[see e.g. \cite{Gille}] It follows directly from the Steinberg relations that
  $\{x,x\}=\{x,-1\}$. The equality $\{x,x\}=\{x,-1\}=0$ happens
  \begin{itemize} \item in $\KM_2(k)/p$ if $p \neq 2$. \item in
    $\KM_2(k)/2$ if $-1 \in (k^\times)^2$. \item in $\KM_2(k)$ if $k$
    has characteristic $2$. \end{itemize}

\end{remark}

Let $K$ be a field equipped with a discrete valuation $v : K^\times \r
\Z$, and $\mathcal{O}_v$ its associated discrete valuation ring. Fix a
local parameter $\pi$ and let $\kappa$ be its residue field. Then, for
each $n \geq 1$, there exists a unique set of homomorphisms
$\partial_v : \KM_n(K) \r \KM_{n-1}(\kappa)$ and $s_\pi: \KM_n(K) \r
\KM_n(\kappa)$ satisfying
$\partial_v(\{\pi,u_2,\ldots,u_n\})=\{\bar{u}_2,\ldots,\bar{u}_n\}$
and
$s_\pi(\{\pi^{i_1}u_1,\ldots,\pi^{i_n}u_n\})=\{\bar{u}_1,\ldots,\bar{u}_n\}$
for all units $u_1,\ldots,u_n$ in $\mathcal{O}_v$, where $\bar{u}_i$
denotes the image of $u_i$ in $\kappa$. The maps $\partial_v$ are
called the residue maps and the maps $s_\pi$ are called the
specialization maps.  The specialization maps depend on the choice of
a local parameter, whereas the residue maps don't. It is easy to see
that these maps induce well-defined maps on the quotients $\KM_*/p$.
Moreover, for any $x \in \KM_*(K)$, they are related by the formula :
$$s_\pi(x)=\partial_v \big(\{-\pi\}\cdot x\big).$$

Let $k(t)$ be the function field over $k$ in one variable. A closed
point $P$ in the projective line $\P_k^1$ over $k$ determines a
discrete valuation on $k(t)$ and can be viewed as an irreducible
polynomial in $k[t]$ (for $P=\infty$, take $P=t^{-1}$).  Let
$\partial_P$ (resp.  $s_P$) be the corresponding residue map (resp.
specialization map) and $\kappa_P$ the residue field corresponding to
the valuation induced by $P$. Then, we have Milnor's exact sequence
$$0 \r \KM_n(k) \r \KM_n(k(t)) \stackrel{\oplus
  \partial_P}{\longrightarrow} \bigoplus_{P\in \P^1_k \backslash
  \{\infty\}} \KM_{n-1}(\kappa_P) \r 0$$ where the injective arrow is
induced by the inclusion of fields $k \subseteq k(t)$. Moreover, this
sequence is split by $s_\infty$. This yields an exact Milnor sequence
mod $p$ for any prime number $p$.

Kummer theory defines a map in Galois cohomology $\partial : k^\times
\r H^1(k,\mu_m)$ where $\mu_m$ is the group of $m^\mathrm{th}$ roots
of unity in a fixed separable closure of $k$. Consider the cup-product
$H^1(k,\mu_m) \otimes \ldots \otimes H^1(k,\mu_m) \r
H^n(k,\mu_m^{\otimes n})$. We get a map $\partial^n : k^\times \otimes
\ldots \otimes k^\times \r H^n(k,\mu_m^{\otimes n})$. Bass and Tate
prove in \cite{BassTate} that this map factors through $\KM_n(k)$ and
yields a map $h_{k,m}^n : \KM_n(k) \r H^n(k,\mu_m^{\otimes n})$. The
map $h_{k,m}^n$ is called the Galois symbol.  The Bloch-Kato
conjecture asserts that the Galois symbol induces an isomorphism
$\KM_n(k)/m \r H^n(k,\mu_m^{\otimes n})$ for all $n\geq 0$, all fields
$k$ and all integer $m$ prime to the characteristic of $k$. The case
when $m$ is a power of $2$ is known as Milnor's conjecture and has
been proven by Voevodsky in \cite{Voevodsky}. The case $n=0$ is
trivial, the case $n=1$ is just Kummer theory and Hilbert 90, and the
case $n=2$ is known as the Merkurjev-Suslin theorem (cf. \cite{MS}).
Rost and Voevodsky have announced a proof of the Bloch-Kato conjecture
in the general case, see \cite{Voevodsky2} and \cite{Weibel}. The
proof relies on the existence of reduced power operations. \medskip

In the sequel, $p$ will always denote a prime number, cup product in
Milnor $K$-theory will be denoted by ``$\cdot$'' and by definition the group
$\KM_i(k)$ will be $0$ as soon as $i<0$. We also write $\KM_*(k)$ for $\bigoplus_{i \geq 0} \KM_i(k)$. By ring, we mean commutative
ring with unit.

\section{Divided powers}

\subsection{Existence of divided powers in Milnor $K$-theory}

In \cite{Kafields}, Kahn mentions the existence of divided powers in
all cases of the following proposition. However, we recall the
construction of divided powers, as it will prove to be useful for the
determination of our operations.

\begin{definition}
  Let $n$ be a non-negative integer and $F$ any field. Let $x =
  \sum_{r=1}^l s_r$ be a sum of $l$ symbols of degree $i$ in some
  Milnor $K$-group. A divided power of $x$ is
$$\gamma_n(x)= \sum_{1 \leq l_1 < \ldots < l_n \leq l} s_{l_1} \cdot \ldots
\cdot s_{l_n}.$$ 
\end{definition}

Of course, it is not clear that  $\gamma_n$ should give a well-defined
map $ \KM_i(k)/p \rightarrow \KM_{ni}(k)/p$. However, we have

\begin{proposition} \label{divpowers} $\gamma_0 =1$ and $\gamma_1 =
  \mathrm{id}$ are always well-defined.

  \begin{enumerate} \item If $i$ is even $\geq 2$ and $p$ is an odd
    prime number, then there exists a divided power $$\gamma_n :
    \KM_{i}(F)/p \r \KM_{ni}(F)/p.$$

  \item If $-1 \in (F^\times)^2$, then for all $i\geq 2$,
    there exists a divided power $$\gamma_n : \KM_{i}(F)/2 \r
    \KM_{ni}(F)/2.$$

  \item If $i$ is even $\geq 2$ and $\mathrm{char} \ F = 2$, then
    there exists a divided power $${\gamma}_n : \KM_{i}(F) \r
    \KM_{ni}(F).$$
  \end{enumerate}

\end{proposition}

\noindent \textbf{Proof.} We are going to prove that ${\gamma}_n$ as
defined explicitly doesn't depend on how we write $x$ as a sum of symbols. We
will give the proof only in the first case. The two remaining cases
can be proven exactly the same way, once one remarks that the
conditions ($i$ even and $p$ odd), ($p=2$ and $-1 \in (F^\times)^2$)
and ($i$ even and $\mathrm{char} \ F = 2$) are here to force :
\begin{itemize} \item $\{a,a\}=0$ \item The algebras $\bigoplus_{n
    \geq 0} \KM_{ni}(F)/p$, $\bigoplus_{n \geq 0} \KM_{ni}(F)/2$ and
  $\bigoplus_{n \geq 0} \KM_{ni}(F)$ are commutative. \end{itemize}

\noindent Let $M_{i,F}$ be the free $\Z$-module generated by elements
in $(F^\times)^i$, and define \begin{center} $\begin{array}{rcl} \
    \tilde{\Gamma}_n : M_{i,F} & \rightarrow & \KM_{ni}(F)/p \\ \
    \sum_{r=1}^l n_r s_r & \mapsto & \sum_{1 \leq l_1 < \ldots < l_n
      \leq l} n_{l_1}s_{l_1} \cdot \ldots \cdot n_{l_n}s_{l_n} \\
  \end{array}$ \end{center}

The map $\tilde{\Gamma}_n$ is well-defined because of the
commutativity of $\bigoplus_{n \geq 0} \KM_{ni}(F)/p$ and we want to
show that $\tilde{\Gamma}_n$ factors through $(F^\times)^{\otimes i}$.
For this purpose, it is enough to show that $\tilde{\Gamma}_n$ takes
the same value on each equivalence class for the quotient map $M_{i,F}
\r (F^\times)^{\otimes i}$. We notice that, given $x$ and $y$ in
$M_{i,F}$, the sum formula $\tilde{\Gamma}_n(x+y) = \sum_{j=0}^n
\tilde{\Gamma}_j(x) \tilde{\Gamma}_{n-j}(y)$ holds. Thus, if $x$ and
$x+y$ have same image in $(F^\times)^{\otimes i}$, we want to show
that $\tilde{\Gamma}_n(x+y) = \sum_{j=0}^n \tilde{\Gamma}_j(x)
\tilde{\Gamma}_{n-j}(y) = \tilde{\Gamma}_n(x)$. For this, it is enough
to prove that $\tilde{\Gamma}_j(y)=0$ for all $j \geq 1$ and $y$
mapping to $0$ in $(F^\times)^{\otimes i}$. Still using the sum
formula, it is enough to prove that for any $1\leq j \leq i$, elements
of the form $(a_1, \ldots, a_{j-1},b,a_{j+1},\ldots,a_i) + (a_1,
\ldots, a_{j-1},c,a_{j+1},\ldots,a_i) - (a_1, \ldots,
a_{j-1},bc,a_{j+1},\ldots,a_i)$ map to $0$ under $\tilde{\Gamma}_n$.
This happens, when $i\geq 2$, because $\{a,a\}=0$ in $\KM_2(F)/p$ for
all $a \in F^\times$.

\noindent Therefore, we get a map $\Gamma_n : (F^\times)^{\otimes i}
\r \KM_{ni}(F)/p$ that satisfies the sum formula.

\noindent Now, $\Gamma_n$ factors through $\tilde{\gamma}_n : \KM_i(F)
\r \KM_{ni}(F)/p$. Let $x=x_1+x_2$ with $x_2$ a pure tensor satisfying
Steinberg's relation. Then $\Gamma_n(z)=\Gamma_n(z_1) +
\Gamma_{n-1}(z_1)\Gamma_1(z_2) + \ldots$, but $\Gamma_n(z_2)$ is
clearly $0$ for $n \geq 1$. Hence a map $\tilde{\gamma}_n : \KM_{i}(F)
\r \KM_{ni}(F)/p$.

\noindent Finally, $\tilde{\gamma}_n$ factors through ${\gamma}_n :
\KM_i(F)/p \r \KM_{ni}(F)/p$ as easily seen. \qed \medskip

\noindent In all cases of Proposition \ref{divpowers}, the divided powers
satisfy the following properties. Moreover, they are the only set of maps
satisfying such properties. See e.g. \cite[Theorem 2]{Kafields}.

\begin{properties} \label{prop} 
  \begin{enumerate}
  \item $\gamma_0(x)=1$, $\gamma_1(x)=x$.
  \item $\gamma_n(xy)=x^n \gamma_n(y)$.
  \item $ \gamma_m(x) \gamma_n(x) = \begin{pmatrix} m+n \\ n
    \end{pmatrix} \gamma_{m+n}(x)$.
  \item $\gamma_n(x+y) = \sum_{i=0}^n \gamma_i(x)\gamma_{n-i}(y)$.
  \item $\gamma_m (\gamma_n(x)) = \frac{(nm)!}{m!n!^m}
    \gamma_{nm}(x)$.
  \item $\gamma_n(s)=0$ if $n\geq 2$ and $s$ is a symbol.
  \end{enumerate}
\end{properties}

All these properties imply that Milnor $K$-theory is a divided power algebra in the sense of Revoy \cite{Revoy}.

\subsection{Divided powers and length}

\begin{definition} The \emph{length} of an element $x \in \KM_i(F)/p$
  (resp. $\KM_i(F)$) is the minimum number of symbols appearing in any
  decomposition of $x$ as a sum of symbols.
\end{definition}

\begin{remark} If $x \in \KM_i(F)/p$ has length $l$ and if the divided
  power $\gamma_n$ is well-defined on $\KM_i(F)/p$, then $n >l$
  implies $\gamma_n(x)=0$. That is, $\gamma_n$ vanishes on elements of
  length strictly less than $n$.
\end{remark}

As was noted by Kahn in \cite{Kafields}, the existence of divided
powers implies

\begin{proposition}

  If there exists an integer $l$ such that the length of any element
  in $\KM_2(F)/p$ (resp. $\KM_2(F)/2$, $\KM_2(F)$) is $\leq l$, then
  for all $n \geq 2l+2$, we have 
  \begin{itemize} 
  \item $\KM_n(F)/p = 0$, if $p$ is odd.
  \item $\KM_n(F)/2 = 0$, if $-1$ is a square in $F^\times$. 
  \item $\KM_n(F)=0$, if $\mathrm{char} \ F =2$. 
  \end{itemize}

\end{proposition}

\noindent \textbf{Proof.} Let $x = \{x_1,\ldots,x_n\}$ be a symbol in
$\KM_n(F)/p$ with $p$ odd and $n \geq 2l+2$. Then $x =
{\gamma}_{l+1}\big(\{{x}_1,{x}_2\} + \ldots +
\{{x}_{2l+1},{x}_{2l+2}\}\big)\cdot \{x_{2l+3},\ldots ,x_n\}$. By
assumption $\{{x}_1,{x}_2\} + \ldots + \{{x}_{2l+1},{x}_{2l+2}\}$ has
length at most $l$ and maps therefore to $0$ under ${\gamma}_{l+1}$.
The two other cases are similar. \qed

\begin{examples} In \cite{Becher2}, Becher shows that if $\KM_1(F)/p =
  F^\times / (F^\times)^p$ is finite of order $p^m$, then the length
  of elements in $\KM_2(F)/p$ is always less or equal than
  $\frac{m}{2}$ if $p$ is odd, and is always less or equal than
  $\frac{m+1}{2}$ if $p=2$. Hence, if $\KM_1(F)/p = F^\times /
  (F^\times)^p$ is finite of order $p^m$, the higher Milnor $K$-groups
  $\KM_n(F)/p$ are zero whenever $n > 2[\frac{m}{2}] +1$ if $p$ is odd
  and $n > 2[\frac{m+1}{2}]+1$ if $p=2$. It is worth saying that
  Becher also shows that these upperbounds are
  sharp. \end{examples}

\subsection{Stiefel-Whitney classes of a quadratic form}

Let $q$ be a quadratic form of rank $r$ over a field $F$ of
characteristic $\neq 2$. Then $q$ admits a diagonal form $\langle
a_1,\ldots,a_r\rangle$.  The total Stiefel-Whitney class of $q$ is
defined to be $w(q) = (1+\{a_1\})\ldots (1+\{a_n\}) \in \KM_*(F)/2$.
In \cite{Milnor}, Milnor shows that $w(q)$ is well-defined and doesn't
depend on a particular choice of a diagonal form for $q$. The
$k^\mathrm{th}$ Stiefel-Whitney class $w_k$ is defined to be the
degree $k$ part of $w$.

\begin{proposition}[Milnor, Becher] \label{SW}

  We have $w_1w_2=w_3$ and more generally if $n = \sum \varepsilon_i
  2^i$ is the binary decomposition of $n$, then $w_n= \prod_{i,
    \varepsilon_i=1} w_{2^i}$. Also, when $-1 \in (F^\times)^2$,
  $w_{2n}=\gamma_n(w_2)$ and so $w_{2n+1} = w_1 \cdot \gamma_n(w_2)$.
  We also have $$w_n = w_1^{\varepsilon_0} \cdot \prod_{i\geq 1,
    \varepsilon_i=1} \gamma_{2^{i-1}}(w_2).$$

\end{proposition}

\noindent \textbf{Proof.} The first point is proved in \cite{Milnor}.
The last point was pointed out by Becher in \cite[paragraph 9]{Becher} and is a
direct consequence of the existence of divided powers mod $2$ when
$-1$ is a square in $F^\times$, and of the explicit formula defining
both the Stiefel-Whitney classes and the divided powers. \qed \medskip

This result confirms that the invariants $w_1$ and $w_2$ of a
quadratic form are important. In the literature, $w_1$ is often
referred to as the determinant, and $w_2$ as the Hasse invariant.
Also, the Witt invariant can be expressed in terms of the determinant
and the Hasse invariant (see e.g.  \cite{Lam}, Proposition V.3.20). A
natural question is to ask whether or not a non-degenerate quadratic
form is determined, up to isometry, by its total Stiefel-Whitney
class. This has been answered by Elman and Lam in \cite{EL1} : let $F$
be a field of characteristic not equal to $2$, with $W(F)$ its Witt
ring of anisotropic quadratic forms, and $IF$ its ideal of
even-dimensional forms. Let $I^nF$ denote the $n$th power $(IF)^n$.
Then, the equivalence class of a non-degenerate quadratic form over
$F$ is determined by its dimension and Stiefel-Whitney invariant if
and only if $I^3F$ is torsion free (as an additive abelian group).
This is for example the case of the field of real numbers. However, a
real non-degenerate quadratic form of given rank is not solely
determined by $w_1$ and $w_2$ and as such, the higher Stiefel-Whitney
classes do carry a little information beyond what $w_1$ and $w_2$
give. On the other hand, proposition \ref{SW} shows that when $-1$ is
a square in the base field $F$, the higher Stiefel-Whitney classes are
completely determined by $w_1$ and $w_2$. This helps to explain why
the classes $w_i$ for $i \geq 3$ have played very little role in
quadratic form theory. Also, in general, it is known that two
non-degenerate quadratic forms $q$ and $q'$ of same dimension $\leq 3$
are isometric if and only if they have same $w_1$ and $w_2$ (see
\cite{Lam}, Proposition V.3.21.). Finally, Elman and Lam gave a
description of fields for which non-degenerate quadratic forms of
given dimension are totally determined, up to isometry, by their
determinant $w_1$ and Hasse invariant $w_2$. This happens if and only
if $I^3F=0$ (cf.  \cite{EL2}).\medskip

By Milnor's conjecture (proven by Voevodsky in \cite{Voevodsky}),
proposition \ref{divpowers} gives divided powers in Galois cohomology.
Let $Et_n$ be the functor that associates to any field $F$ over $k_0$
the set of \'etale algebras of rank $n$ over $F$. In \cite{Galois}, it
is proven that the $H^*(k_0,\Z/2)$-module
$\mathrm{Inv}_{k_0}(Et_n,\Z/2)$ of natural transformations from the
functor $\mathrm{Et}_n$ over $k_0$ to the functor $H^*(-,\Z/2)$ over
$k_0$ is free with basis the Galois-Stiefel-Whitney classes $1,
w_1^{gal}, \ldots, w_m^{gal}$, $m=[\frac{n}{2}]$. Moreover,
$w_i^{gal}=0$ for $i>m$. Given an \'etale algebra $E$ over $k_0$, we
can consider the non-degenerate quadratic form $q_E$ on $E$ viewed as
a $k_0$-vector space, defined as $q_E(x) = \mathrm{Tr}_{E/k_0}(x^2)$.
Therefore, we have some invariants, coming from the Stiefel-Whitney
classes of quadratic forms. In \cite{KahnSW}, Kahn proves that for $E
\in Et_n(k_0)$, $w_i^{gal}(E) = w_i(q_E)$ if $i$ is odd and
$w_i^{gal}(E) = w_i(q_E) + (2)\cdot w_{i-1}(q_E)$ if $i$ is even. When
$-1$ is a square in $k_0$, we get that the higher
Galois-Stiefel-Whitney invariants are determined by $w_1^{gal}$ and
$w_2^{gal}$ in the following way : $w_{2n+1}^{gal} = w_1^{gal} \cdot
\gamma_n(w_2^{gal}-\{2\} \cdot w_1^{gal})$ and
$w_{2n}^{gal}=\gamma_n(w_2^{gal}-\{2\} \cdot w_1^{gal}) + \{2\} \cdot
w_1^{gal} \cdot \gamma_{n-1}(w_2^{gal}-\{2\} \cdot
w_1^{gal})$.\medskip

However, there are some examples where the divided powers act
trivially.  For instance, MacDonald computes, for $n$ odd $\geq 3$,
the mod $2$ cohomological invariants for the groups $SO(n)$, $\Z/2
\ltimes PGL(n)$, $PSp(2n)$, and $F_4$. These correspond, for
$r=0,1,2$, and $3$ respectively, to $\mathrm{Inv}_{k_0}(J_n^r,\Z/2)$,
the group of mod $2$ invariants for odd degree $n\geq 3$ Jordan
algebras with associated composition algebra of dimension $2^r$
($0\leq r \leq 3$) over the base field $k_0$ of characteristic
supposed to be different from $2$. Such algebras are known to be of
the form $H(C,q) = \{x \in M_n(C), \ B_q^{-1}\bar{x}^t B_q = x \}$ for
$q$ an $n$-dimensional quadratic form of determinant $1$ with
associated bilinear form $B_q$ and $C$ a composition algebra over
$k_0$ of dimension $2^r$. The composition algebra comes with a norm
form $\varphi$, which turns out to be a Pfister form. The group of
invariants for $r$-Pfister forms ($r>0$) is the free
$H^*(k_0,\Z/2)$-module generated by $1$ and $e_r$, where
$e_r(\langle\langle a_1, \ldots, a_r \rangle\rangle) =
(a_1)\cdot\ldots \cdot (a_r)$. Write $J= \varphi \otimes q$ for a
Jordan algebra $J$, then we have invariants $v_i = e_r \otimes w_{2i}$
and it is shown in \cite{MacDonald} that
$\mathrm{Inv}_{k_0}(J_n^r,\Z/2)$ is the free $H^*(k_0,\Z/2)$-module
generated by $1, v_0, \ldots, v_m$, with $m$ satisfying $n=2m+1$. For
$r>0$, when $-1$ is a square in $k_0$ and because $e_r \cdot e_r$ is
zero, we see that the divided powers vanish on the $v_i$'s.

\subsection{Divided powers in Milnor $K$-theory mod $2$}

In this section, we no longer assume $-1 \in (F^\times)^2$. We define
the map $$\tau_i : \KM_*(F)/2 \r \KM_*(F)/2, \ \ x \mapsto
\{-1\}^{i-1} \cdot x.$$ Let's say a few words about this map.  If $F$
has characteristic $p>0$, then $\{-1\}$ is zero in $\KM_1(F)/2$ if
$p\equiv 1$ mod $4$, and in any case $\{-1,-1\}=0$ since the groups
$\KM_n(\mathbf{F}_q)$ vanish for finite fields $\mathbf{F}_q$ and $n
\geq 2$. So, considering a function field over a finite field, we see
that the maps $\tau_i$ for $i \geq 2$ are neither injective nor
surjective, even when $-1$ is not a square in $F$. If $F$ is a number
field (or a global field), let $r_1$ be the number of real places of
$F$ and denote them by $\sigma_i : F \r \mathbf{R}$. Bass and Tate
show in \cite{BassTate} that for $n\geq 3$, the embeddings $\sigma_i :
F \r \mathbf{R}$ corresponding to the real places of $F$ induce an
isomorphism $\KM_n(F) \stackrel{\oplus \sigma_i}{\longrightarrow}
\bigoplus_{i=1}^{r_1} \KM_n(\mathbf{R})/2 \cong \big( \Z/2
\big)^{r_1}$.  Then, clearly $\KM_n(F)/2 \cong \big( \Z/2 \big)^{r_1}$
for $n \geq 3$. Also, $\KM_1(F)/2$ is countably infinite. This shows
that $\tau_i$ cannot be injective. 

Hence $\ker (\tau_i)$, or equivalently the annihilator ideal of
$\{-1\}^{i-1}$ in $\KM_*(F)/2$, is non-trivial in general.

\begin{proposition} \label{mod2} Let $n$ be an integer $\geq 2$ and
  $F$ any field. Let $y_n$ be in the kernel of $\tau_i$. Then, if
  $s_1, \ldots,s_l$ are symbols in $\KM_i(F)/2$, $$(y_n \cdot
  \gamma_n) (s_1 +\ldots + s_l) = y_n \cdot \sum_{1\leq s_{l_1} <
    \ldots < s_{l_n} \leq l} s_{l_1} \cdot \ldots \cdot s_{l_n}$$ is a
  well-defined map over $\KM_i(F)/2$.
\end{proposition}

\noindent \textbf{Proof.} We proceed exactly the same way as in
Proposition \ref{divpowers}, from which we take up the notations. The
map $y_n \cdot \tilde{\Gamma}_n : M_{i,F} \rightarrow \KM_*(F)/2, \
\sum_{r=1}^l n_r s_r \mapsto y_n \cdot \sum_{1 \leq l_1 < \ldots < l_n
  \leq l} n_{l_1}s_{l_1} \cdot \ldots \cdot n_{l_n}s_{l_n}$ is
well-defined due to the commutativity of the $\mathbf{F}_2$-algebra
$\KM_*(F)/2$. As before, $\tilde{\Gamma}_n$ satisfies a sum formula
which we write $y_n \cdot \tilde{\Gamma}_n(x+y) = \sum_{j=0}^n
\tilde{\Gamma}_j(x) \cdot y_n \cdot \tilde{\Gamma}_{n-j}(y)$ for all
$x$ and $y$ in $M_{i,F}$. To prove that $\tilde{\Gamma}_n$ factors
through $(F^\times)^{\otimes i}$, it is enough to show that elements
of the form $y = (a_1, \ldots, a_{i-1},b) + (a_1, \ldots, a_{i-1},c) -
(a_1, \ldots, a_{i-1},bc)$ map to zero under $y_n \cdot
\tilde{\Gamma}_j$ for all $j \geq 1$. This is clear for $j=1$ and
$j>3$. In the case $j=2$, we have $y_n \cdot \tilde{\Gamma}_2 (y) =
y_n \cdot \big( \{a_1, \ldots, a_{i-1},b,a_1, \ldots, a_{i-1},c\} +
\{a_1, \ldots, a_{i-1},b,a_1, \ldots, a_{i-1},bc\} + \{a_1, \ldots,
a_{i-1},c,a_1, \ldots, a_{i-1},bc\} \big)$.  Notice that $$\{a_1,
\ldots, a_{i-1},a_1, \ldots, a_{i-1}\} = \{-1\}^{i-1} \{a_1, \ldots,
a_{i-1}\}$$ to conclude $y_n \cdot \tilde{\Gamma}_2 (y) =0$.  The case
$j=3$ is similar.

\noindent Thus, we have a well-defined map $y_n \cdot \Gamma_n :
(F^\times)^{\otimes i} \r \KM_*(F)/2$, and as in the proof of
Proposition \ref{divpowers}, it factors through a well-defined map
$y_n \cdot \gamma_n : \KM_i(F)/2 \r \KM_*(F)/2$.  \qed \medskip

\begin{example} Consider again quadratic forms over $F$, but this time
without assuming $-1$ is a square in $F$, and their Stiefel-Whitney
invariants. If $y \in \KM_*(F)/2$ is such that $\{-1\}\cdot y = 0$, then
$$y\cdot w_{2n} = (y \cdot \gamma_n)(w_2).$$ If we write $w_2$ as an
ordered sum of symbols $\sum_i s_i$, we can consider $\gamma_n(w_2)$
to be $\sum_{i<j} s_i \cdot s_j$. Of course, this may not be
independent on the choice of the $s_i$'s.  However, we have for all $y
\in \kr \tau_2$, $y\cdot (\gamma_n(w_2) - w_{2n})=0$. Hence,
$\gamma_n(w_2)-w_{2n}$ must be in the subgroup $G$ of $\KM_{2n}(F)/2$
consisting of elements $z$ such that $z\cdot y=0$ for all $y\in \kr
\tau_2$. In particular, $G$ contains $\{-1\}\cdot \KM_{2n-1}(F)/2$.
This means that knowing $w_1$ and $w_2$ gives some restriction on the
possible values of the higher Stiefel-Whitney classes even when $-1$
is not a square.

\end{example}

\subsection{Divided powers in integral Milnor $K$-theory}

In this section, $\tau_i$ is the map on integral Milnor $K$-theory
$\KM_*(F) \r \KM_*(F), \ x \mapsto \{-1\}^{i-1} \cdot x$. The same
examples as in the previous section show that this map is not
necessarily injective nor surjective.

\begin{proposition} \label{integraleven} Let $n$ and $i$ be integers
  $\geq 2$ with $i$ even, and $F$ any field. Let $y_n$ be an element
  in the kernel of $\tau_i$. Then, if $s_1, \ldots, s_l$ are symbols
  in $\KM_i(F)$, $$(y_n \cdot \gamma_n) (s_1 +\ldots + s_l) = y_n
  \cdot \sum_{1\leq s_{l_1} < \ldots < s_{l_n} \leq l} s_{l_1} \cdot
  \ldots \cdot s_{l_n}$$ is a well-defined map over $\KM_i(F)$.
\end{proposition} 

\noindent \textbf{Proof.} Same as for Proposition \ref{mod2} since the
algebra $\bigoplus_{r\geq 0} \KM_{2r}(F)$ is commutative. \qed

\begin{proposition} \label{integralodd} Let $n$ and $i$ be integers
  $\geq 2$ with $i$ odd, and $F$ any field. Let $y_n$ be an element in
  the kernel of $\tau_i$, which is of $2$-torsion. Then, if $s_1,
  \ldots,s_l$ are symbols in $\KM_i(F)$, $$(y_n \cdot \gamma_n) (s_1
  +\ldots + s_l) = y_n \cdot \sum_{1\leq s_{l_1} < \ldots < s_{l_n}
    \leq l} s_{l_1} \cdot \ldots \cdot s_{l_n}$$ is a well-defined map
  $\KM_i(F)$ to $\KM_{ni}(F)$.
\end{proposition}

\noindent \textbf{Proof.} Notice that the map $y_n \cdot
\tilde{\Gamma}_n : M_{i,F} \rightarrow \KM_*(F), \ \sum_{l=1}^r n_l
s_l \mapsto y_n \cdot \sum_{1 \leq l_1 < \ldots < l_n \leq r}
n_{l_1}s_{l_1} \cdot \ldots \cdot n_{l_n}s_{l_n}$ is well-defined
because $y_n$ is of $2$-torsion. Now, the proof is the same as for
Proposition \ref{mod2}. \qed

\section{Operations in Milnor $K$-theory of a field}

We start this section with a result that will be of constant use.

\begin{proposition} \label{inj}

  Let $a$ be in $\KM_n(k_0)$ and suppose that for all extension
  $k/k_0$ and for all $x \in \KM_{>0}(k)$ we have $a \cdot x = 0$, then
  $a$ is necessarily $0$ in $\KM_n(k_0)$. Moreover the same result
  holds mod $p$.

\end{proposition}

\noindent \textbf{Proof.} Let $a$ be as in the proposition. Consider
the map $\KM_n(k_0) \r \KM_{n+1}(k_0(t)), \ a \mapsto \{t\} \cdot a$.
This map is injective since it admits a left inverse, namely the
residue map $\partial_0^\mathrm{M}$ at the point $0 \in \P^1_{k_0}$.
Indeed, we have the formula $\partial_0^\mathrm{M}(\{t\}\cdot a)=a$.
The residue map, as defined in \cite[Chapter $7$]{Gille}, is a
homomorphism and hence induces a well-defined residue map mod $p$,
$\KM_{n+1}(k_0(t))/p \r \KM_{n}(k_0)/p$. Hence, the same arguments
apply in the mod $p$ case. \qed

\begin{definition}
  An operation $\varphi : \KM_i/p \r \KM_*/p$ over a field $k_0$ is a
  natural transformation from the functor $ \KM_i/p :
  \mathsf{Fields}_{/k_0} \r \mathsf{Sets}$ to the functor $ \KM_*/p :
  \mathsf{Fields}_{/k_0} \r \mathbf{F}_p - \mathsf{Algebras}$. In
  other words, it is a set of functions $\varphi : \KM_i(k)/p \r
  \KM_*(k)/p$ defined for all extensions $k$ of $k_0$ such that for
  any extension $l$ of $k$, the following diagram commutes:
  \begin{center} $ \xymatrix{\KM_i(l)/p \ar[r]^\varphi & \KM_*(l)/p \\
      \KM_i(k)/p \ar[r]^{\varphi} \ar[u] & \KM_*(k)/p \ar[u] }$
  \end{center}

\end{definition}

\begin{ex}
  Divided powers are indeed operations in the above sense
  (when they are well-defined). So, if for instance $p$ is odd and $i$
  is even, any sum of divided power operations with coefficients in
  $\KM_*(k_0)/p$ gives an operation $\KM_i/p \r \KM_*/p$ over $k_0$. Our
  main theorems say that this gives all the possible operations.
\end{ex}

\begin{ex}
  Suppose $p$ is odd, $i\geq 2$ is even and $k$ is an extension of
  $k_0$. The map $\KM_i(k)/p \r \KM_{2i}(k)/p, \ x \mapsto x^2$
  defines an operation over $k_0$. It is easy to check that this
  operation corresponds to $2\cdot \gamma_2$. More generally, it is
  straightforward to check that any map of the form $x \mapsto x^q$
  defines an operation $\KM_i(k)/p \r \KM_{qi}(k)/p$ and that it is a
  sum of divided powers. Of course, this is a particular case of
  Theorem \ref{main}. More precisely, $x^q$ is equal to $0$ if $i$ is
  odd and is equal to $q!\gamma_q(x)$ if $i$ is even.
\end{ex}

\begin{definition}\label{spedef}
  Let $k_0$ be any field and $K$ an extension of $k_0$ endowed with a
  discrete valuation $v$ such that its valuation ring $R = \{x \in K,
  v(x) \geq 0\}$ contains $k_0$, so that the residue field $\kappa$ is
  an extension of $k_0$. We say that specialization maps commute with an
  operation $\varphi : \KM_i \longrightarrow \KM_*$ over $k_0$ if for
  any extension $K/k_0$ as above, we have a commutative diagram
  \begin{center} $ \xymatrix{\KM_i(K)/p \ar[r]^{
        \ \varphi} \ar[d]_{s_\pi} & \KM_*(K)/p \ar[d]_{s_\pi} \\
      \KM_i(\kappa)/p \ar[r]^{ \ \varphi} & \KM_*(\kappa)/p }$
  \end{center} where $\pi$ is any uniformizer for the valuation $v$.
\end{definition}

\begin{ex}
  Divided power operations over $k_0$ do commute with specialization
  maps. This is clear from the definition of specialization maps.
\end{ex}

\subsection{Operations $\KM_1/p \times \ldots \times \KM_1/p \r \KM_*/p$}

The following theorem is essential in the determination of operations
$\KM_i/p \r \KM_*/p$ for $i\geq 2$.

\begin{theorem}  \label{theorem}

  Let $k_0$ be any field and let $p$ be a prime number. The algebra of
  operations
$$ \underbrace{ \KM_1(k)/p \times \ldots \times \KM_1(k)/p}_{r \
  \mathrm{times}} \longrightarrow \KM_*(k)/p$$ for fields $k \supseteq
k_0$ and commuting with field extension over $k_0$ is the free module
over $\KM_*(k_0)/p$ with basis the operations $\big( \{a_1\}, \ldots,
\{a_r\}\big) \mapsto \{a_{i_1},\ldots,a_{i_s}\}$ for all subsets $1
\leq i_1 < \ldots < i_s \leq r$.

\end{theorem} 

For example, given an operation $\psi : \KM_1(k)/p \times \KM_1(k)/p
\r \KM_5(k)/p$, there exist $a \in \KM_5(k_0)/p$, $b_1$ and $b_2 \in
\KM_4(k_0)/p$ and $c \in \KM_3(k_0)/p$, such that for any
$(\{x\},\{y\}) \in \KM_1(k)/p \times \KM_1(k)/p$, we have $\psi
(\{x\},\{y\}) = a + b_1 \cdot \{x\} + b_2 \cdot \{y\} + c \cdot
\{x,y\}$.\medskip


\noindent \textbf{Proof of Theorem \ref{theorem}.} The proof goes in
three steps. In the first step, we show that an operation $\KM_1/p \r
\KM_*/p$ over $k_0$ is determined by the image of $\{t\} \in
\KM_1(k_0(t))/p$ where $t$ is a transcendental element over $k_0$. In
the second step, we determine the image of $\{t\}$. Finally, in the
last step we conclude by induction on the number $r$ of factors. Let
$\varphi : \KM_1/p \r \KM_*/p$ be an operation over $k_0$. \medskip

\emph{Step 1.} The operation $\varphi : \KM_1/p \r \KM_*/p$ over $k_0$
is determined by the image $\varphi(\{t\})$ of $\{t\} \in
\KM_1(k_0(t))/p$ in $\KM_*(k_0(t))/p$, for $t$ transcendental over
$k_0$. Indeed, consider a field extension $k/k_0$ and an element $e
\in k$. If $e$ is not algebraic over $k_0$, then $\varphi(\{e\})$ is
the image in $\KM_*(k)/p$ of the element $\{e\} \in \KM_1(k_0(e))/p$.
If $e$ is algebraic and if $k$ possesses a transcendental element $t$
over $k_0$, then $et^p$ is transcendental over $k_0$. Also, in
$\KM_1(k)/p$, $\{et^p\} = \{e\}$, and so $\varphi(\{e\})$ is
determined by $\varphi(\{et^p\})$. Finally, if $e \in k$ is algebraic
over $k_0$, consider the function field $k(t)$ and the commutative
diagram \begin{center} $ \xymatrix{\KM_1(k)/p \ar[d]_\varphi \ar[r]^{i
      \ \ } & \KM_1(k(t))/p \ar[d]_\varphi \\ \KM_*(k)/p \ar[r]^{i \ \
    } & \KM_*(k(t))/p }$ \end{center} where $i$ is the map induced by
the inclusion of fields $k \subset k(t)$. We can write
$\varphi(i(\{e\}))=i(\varphi(\{e\}))$. The element $\varphi(i(\{e\}))$
is determined by the previous case. By Milnor's exact sequence, $i$ is
injective. Therefore, $\varphi(\{e\})$ is uniquely determined. \medskip

\emph{Step 2.} The element $\varphi(\{t\}) \in \KM_*(k_0(t))/p$ has
residue $0$ for all residue maps corresponding to closed points in
$\P_{k_0}^1 \backslash \{0,\infty\}$. 

To prove this, let $X$ be a transcendental element over $k_0(t)$ and denote by $\iota$ the homomorphism in Milnor $K$-theory induced by any inclusion of field $k \subset k(X)$ for $k$ any extension of $k_0$. By definition, the map $\iota$ commutes with $\varphi$. We start by proving that $\iota \circ \varphi(\{t\}) \in \KM_*(k_0(t,X))/p$ has only residue at polynomials with coefficients in $k_0$. 

Recall that Milnor's exact sequence $$0 \r \KM_n(k_0) \r \KM_n(k_0(t)) \stackrel{\oplus
  \partial_P}{\longrightarrow} \bigoplus_{P\in \P^1_k \backslash
  \{\infty\}} \KM_{n-1}(\kappa_P) \r 0$$ is split. Write $\psi_P$ for a splitting map to $\partial_P$ so that, for any $x \in \KM_*(k_0(t))$, we have $$x =  s_\infty(x) + \sum_{P\in \P^1_k \backslash
  \{\infty\}} \psi_P \circ \partial_P(x).$$  
For $P$ a closed point in $\P_{k_0}^1 \backslash \{\infty\}$, it is possible to  view it as a monic
non-constant irreducible polynomial in $k_0[t]$. Let's write $\partial_{P\otimes k_0(X)}$ (resp. $\psi_{P\otimes k_0(X)}$) for the residue map (resp. a splitting to the residue map) at the
polynomial $P \in k_0[t]$ seen as a polynomial in $k_0(X)[t]$ via the
obvious inclusion of fields $k_0 \subset k_0(X)$. Then, we have the following commutative diagrams (See Lemma
\ref{com}).
 \begin{center} $ \xymatrix{\KM_*(k(t)) \ar[d]_{\partial_{P}}
      \ar[r]^{\iota \ } & \KM_*(k(t,X)) \ar[d]^{\partial_{P \otimes k_0(X)}} \\
      \KM_{*-1}(\kappa_P) \ar[r]^{\iota \ } &
      \KM_{*-1}(\kappa_P(X))}$ $\ \ \ \ \ \ \ \ \ $ $ \xymatrix{\KM_{*-1}(k(t))
      \ar[d]_{\psi_P} \ar[r]^{\iota \ } & \KM_{*-1}(k(t,X))
      \ar[d]^{\psi_{P\otimes k_0(X)}} \\ \KM_{*}(\kappa_P) \ar[r]^{\iota \ } &
      \KM_{*}(\kappa_P(X)) }$

   \end{center}
Therefore, $$\iota \circ \varphi(\{t\}) =  s_\infty \big( \iota \circ \varphi(\{t\})\big) + \sum_{P\in \P^1_k \backslash
  \{\infty\}} \psi_{P \otimes k_0(X)} \circ \partial_{P \otimes k_0(X)} \big( \iota \circ \varphi(\{t\})\big),$$ which shows our claim.



Now, given a polynomial $Q$ in $k_0[t]$, let's write $Q_X$ for the
polynomial in $k_0(X)[t]$ defined by $Q_X(t)=Q(tX^p)$. Then, via the
isomorphism $k_0(t) \stackrel{\simeq}{\rightarrow} k_0(tX^p), \ t
\mapsto tX^p$, we get that $\varphi(\{tX^p\})$ has non-zero residues
only at polynomials of the form $Q_X$.

Exploiting the fact that $\{t\}=\{tX^p\}$ in $ \KM_1(k_0(t,X))/p$, we
deduce that if $\varphi(\{t\})=\varphi(\{tX^p\})$ has non-zero residue
at some polynomial $Q_X$ as above, then $Q_X$ must come from a
polynomial $P \in k_0[t]$. Concretely, we must have $$Q(tX^p) = \alpha
(X) P(t), \ \ \mathrm{for} \ \mathrm{some} \ \alpha(X) \in k_0(X), \ P
\in k_0[t].$$ Having in mind that $P$ is monic irreducible, this
implies $P=t$ as easily seen.

So, we have proven that $\varphi(\{t\})$ is unramified outside
$\{0,\infty\}$. Writing $a$ for $\partial_0 \varphi(\{t\})$,
$\varphi(\{t\}) - a\cdot \{t\}$ is then unramified on $\P^1 \backslash
\{\infty\}$. By Milnor's exact sequence, this implies that this
element comes from an element $b \in \KM_*(k_0)/p$. Therefore, we have
$\varphi(\{t\}) = a\cdot \{t\} +b$. Combining with step 1, this tells
us that there exist $a$ and $b$ in $\KM_*(k_0)/p$, such that for any
field extension $k/k_0$ and for any $x \in \KM_1(k)/p$,
$\varphi(x)=a\cdot x + b$. This defines an operation as one can easily
check. \medskip

\emph{Step 3.} We conclude by induction on $r$. The case $r=1$ has
been treated in steps 1 and 2. Now, write $\big(\KM_1(k)/p \big)^r =
\big(\KM_1(k)/p \big)^{r-1} \times \KM_1(k)/p$ and write $(x,\{a_r\})$ for
the element $(\{a_1\},\ldots,\{a_r\}) \in \big(\KM_1(k)/p \big)^r$. Fix a
field $k/k_0$ and $x \in \big(\KM_1(k)/p \big)^{r-1}$, then
$\varphi(x,\{a_r\})$ defines an operation $\KM_1/p \r \KM_*/p$ over $k$.
Steps 1 and 2 yield the existence of elements $c_x$ and $d_x$ in
$\KM_*(k)/p$ such that $\varphi(x,\{a_r\})= c_x\cdot \{a_r\} + d_x$. The maps
$x \mapsto c_x$ and $x \mapsto d_x$ define operations $\big(\KM_1/p
\big)^{r-1} \r \KM_*/p$ over $k_0$. By the induction hypothesis, they
are of the form stated in the theorem. It is now easy to see that
$\varphi(x,\{a_r\})= c_x\cdot a_r + d_x$ is of the required form. It
remains to prove that the operations $(\{a_1\},\ldots,\{a_r\}) \mapsto
\{a_{i_1},\ldots,a_{i_s}\}$ for subsets $1 \leq i_1 < \ldots < i_s
\leq r$ form a free basis. Let $\varphi$ be an operation $\big(\KM_1/p
\big)^{r} \r \KM_*/p$ over $k_0$. By the above, we know that there
exist elements $\lambda_{i_1,\ldots,i_s} \in \KM_*(k_0)/p$ such that
for all field extension $k/k_0$ and all $r$-tuple $(\{a_1\},\ldots,\{a_r\})
\in \big(\KM_1(k)/p \big)^{r}$, $$\varphi(\{a_1\},\ldots,\{a_r\}) =
\sum_{1\leq i_1 < \ldots < i_s \leq r} \lambda_{i_1,\ldots,i_s} \cdot
\{a_{i_1},\ldots,a_{i_s}\}.$$ Assume $\varphi$ is $0$. Fix a subset $1
\leq i_1 < \ldots < i_s \leq r$ and consider the field
$k=k_0(t_{i_1},\ldots,t_{i_s})$ where $t_{i_1},\ldots,t_{i_s}$ are
indeterminates. Let $a$ be the element in $\big(\KM_1(k)/p \big)^r$
with entry $\{t_{i_q}\}$ in the $i_q^\mathrm{th}$ coordinate for all
$q$ and zero elsewhere. Consider also the residue maps corresponding
to the local parameters $t_{i_q}$, $\partial_q :
\KM_*\big(k(t_{i_1},\ldots,t_{i_q})\big)/p \r
\KM_{*-1}\big(k(t_{i_1},\ldots,t_{i_{q-1}})\big)/p$. Then,
$\lambda_{i_1,\ldots,i_s} =
\partial_s \circ \ldots \circ \partial_1( \varphi(x)) = 0$. \qed \medskip

As a consequence of the very nice form of the operations $\KM_1/p \r
\KM_*/p$ over $k_0$, we get

\begin{corollary}
  Let $k_0$ be any field. Then, specialization maps commute with
  operations $\varphi : \KM_1/p \longrightarrow \KM_*/p$ over $k_0$.
\end{corollary} 

\noindent \textbf{Proof.} For $K$ and $v$ as in Definition
\ref{spedef}, specialization maps $s_\pi$ are $\KM_*(k_0)/p$-linear
for any choice of uniformizer $\pi$, as easily seen from their
definition. \qed \medskip

\begin{remark} In \cite{Galois}, Theorem \ref{theorem} is proven in
  Galois cohomology for base fields $k_0$ of characteristic different
  from $p$. The proof relies on the fact that it is possible to show
  first that operations $H^1(-,\Z/p) \r H^*(-,\Z/p)$ over $k_0$
  commute with specialization maps in the above sense. Roughly, this
  is done by proving that the specialization maps admit right inverses
  that are induced by some inclusion of fields. Here, we first
  determine all the operations and obtain \textit{a posteriori} that
  the operations commute with specialization maps. Also, using the
  Faddeev exact sequence for Galois cohomology with finite
  coefficients (\cite[Cor. 6.9.3]{Gille}) and thanks to Kummer theory,
  the proof of Theorem \ref{theorem} translates \emph{mutatis
    mutandis} to the case of operations $H^1(-,\Z/p) \r H^*(-,\Z/p)$
  over $k_0$. Thus, this gives a new proof of \cite[Theorem
  16.4]{Galois}, without assuming the Bloch-Kato conjecture. 
\end{remark}

\begin{remark}
In the case $p=\mathrm{char} \ k_0$ the differential symbols $\psi_k^n
: \KM_n(k)/p \r \nu(n)_k$, $\{x_1,\ldots,x_n\} \mapsto (dx_1/x_1)
\wedge \ldots \wedge (dx_n/x_n)$ are isomorphisms for all $n\geq 0$
(Bloch-Gabber-Kato Theorem) and hence give us the operations for
logarithmic differentials. For more details on the differential
symbols, we refer to \cite[Ch. 9]{Gille}.\medskip
\end{remark}

\subsection{Operations in Milnor $K$-theory mod $p$}

In this section, we determine the group of operations $\KM_i/p \r
\KM_*/p$ over any field $k_0$, except in the case when $p=2$ and $-1$
is not a square in $k_0$. The proof of this last case is postponed to
the next section. Thus we prove here Theorem \ref{main}, and Theorem
\ref{main2} in the particular case when $-1$ is a square in $k_0$.
Proposition \ref{main1} covers the cases when divided powers are
well-defined on $\KM_i(k_0)/p$ (as in Proposition \ref{divpowers}) and
Proposition \ref{odd} deals with the case where $p$ and $i$ are both
odd.

\begin{proposition} \label{main1}

  Let $k_0$ be a field and $p$ a prime number. The algebra of
  operations $\KM_i(k)/p \r \KM_*(k)/p$ commuting with field
  extensions over $k_0$ is

\begin{itemize}

\item If $i=0$, the free $\KM_*(k_0)/p$-module of rank $p$ of
  functions $\mathbf{F}_p \r \KM_*(k_0)/p$. 
\item If $i=1$, the free $\KM_*(k_0)/p$-module of rank $2$, generated
  by $\gamma_0$ and $\gamma_1$.
\item If $i$ is even $\geq 2$ and $p$ is an odd prime, the free
  $\KM_*(k_0)/p$-module generated by the divided powers $\gamma_n$ for
  $n\geq 0$ and the action of one of its element $(y_0,y_1, \ldots)$
  is given by $$x \mapsto y_0 + y_1 \cdot x + y_2 \cdot \gamma_2(x) +
  \ldots + y_l \cdot \gamma_l(x) + \ldots$$
\item If $i\geq 2$, $p=2$ and $-1$ is a square in $k_0^\times$,
  likewise, the free $\KM_*(k_0)/2$-module generated by the divided
  powers $\gamma_n$ for $n\geq 0$.

\end{itemize}
\end{proposition}

\begin{remark} Assume $k_0$ is a field of characteristic $\neq p$ and
  $k$ is any extension of $k_0$. Via the Galois symbol $\KM_i(k)/p \r
  H^i(k,\Z/p(i))$ and under the assumption of the Bloch-Kato
  conjecture, this gives all operations $H^i(-,\Z/p(i)) \r
  H^*(-,\Z/p(*))$ over $k_0$ (by $H^i(k,\Z/p(i))$, we mean the \'etale
  cohomology of the field $k$ with values in $\mu_p^{\otimes i}$).
  Since $k_0$ has characteristic not $p$, we get by Galois descent
  that $H^i(k,\Z/p(j)) = H^i(k,\Z/p(i))\otimes \mu_p^{\otimes (j-i)}$.
  Thus, if $p$ is an odd prime, $i \geq 2$ is even and if $j$ is any
  integer, there is a well-defined divided power operation $\gamma_n :
  H^i\big(-,\Z/p(j)\big) \r H^{ni}\big(-,\Z/p(ni+j-i)\big)$ over
  $k_0$.
\end{remark}

\noindent \textbf{Proof.} In the first case ($i=0$), the functor
$\KM_0/p$ is just the constant functor with value $\mathbf{F}_p$,
hence the result.  The second case is Theorem \ref{theorem}. We now
restrict our attention to the two last cases, that is either $i$ even
and $p$ odd, or $p=2$ and $k_0$ has a square-root of $-1$. By
Proposition \ref{divpowers} the given maps define operations.
Therefore, the inclusion ``$\supseteq$'' holds in each case if we can
prove that the algebra of such operations is a free module over
$\KM_*(k_0)/p$. Let $\psi = y_{l} \cdot \gamma_l + y_{i_1}\cdot
\gamma_{i_1} + \ldots + y_{i_r} \cdot \gamma_{i_r}$ with $l<i_1<\ldots
<i_r$ be an operation on $\KM_i(k)/p$ which is zero, with minimal $l$
such that $y_l \neq 0$. We are going to see that this is
contradictory. If $l=0$, $\psi(0) = y_0 = 0$. If $l\geq 1$, let's
consider the field extension $F=k_0(x_{j,k})_{1\leq j \leq i, 1\leq k
  \leq l}$ in $il$ indeterminates over $k_0$ and the element $x=
\sum_{k=1}^l \{x_{1,k}, \ldots, x_{i,k}\}$ of length at most $l$ in
$\KM_i(F)/p$, then $\psi(x)=y_l \cdot \gamma_l(x)= y_l \cdot
\{x_{1,1}, \ldots, x_{i,l}\}$, which is zero only if $y_l=0$ by
Proposition \ref{inj}.

It is thus enough to prove that given an operation $\varphi : \KM_i/p
\r \KM_j/p$ over $k_0$, $\varphi$ must be of the form $\varphi = y_0 +
y_1 \cdot \mathrm{id} + y_2 \cdot \gamma_2 + \ldots + y_l \cdot
\gamma_l + \ldots$.

Let $k$ be an extension of $k_0$ and let $e$ be an element of
$\KM_i(k)/p$, say of length $\leq l$.  Write $e= \sum_{k=1}^l
\{e_{1,k}, \ldots, e_{i,k}\} \in \KM_i(k)/p$ and adjoin $il$
indeterminates $X_{m,k}$, $1\leq m \leq i, \ 1 \leq k \leq l$, to the
field $k$. Then the same arguments as in Step 1 of the proof of
Theorem \ref{theorem} show that $\varphi(e)$ is determined by
$\varphi\big(\sum_{k=1}^l \{e_{1,k}X_{1,k}^p, \ldots,
e_{i,k}X_{i,k}^p\}\big) \in \KM_i\big(k(X_{m,k})_{1\leq m \leq i, \ 1
  \leq k \leq l}\big)/p$. Moreover, the elements $e_{m,k}X_{m,k}^p$
are independent transcendental elements over the field $k$.

Therefore, given the field $F=k_0(x_{j,k})_{1\leq j \leq i, 1\leq k
  \leq l}$ in $il$ indeterminates over $k_0$, it is enough to study
the image under $\varphi$ of the element $x= \sum_{k=1}^l \{x_{1,k},
\ldots, x_{i,k}\} \in \KM_i(F)/p$ of length $\leq l$. For this
purpose, define the set $E_{i,l} = \{({m,k}), \ 1\leq m \leq i, \ 1
\leq k \leq l \}$ and equip it with the lexicographic order : $(m,k)
\leq (m',k')$ if either $k<k'$ or $k=k'$ and $m\leq m'$. Consider an
element $a$ in $\mathcal{P}(E_{i,l})$ the set of subsets of $E_{i,l}$.
We will write $a_x = \prod_{(m,k)\in a} \{x_{m,k}\}$ for the ordered
product of the $\{x_{m,k}\}$'s for $(m,k) \in a$, and if $a$ is the
empty set, $a_x=1$. With our notations, Theorem \ref{theorem} tells us
that there are unique elements $c_a \in \KM_{j-\#(a)}(k_0)/p$ for
$a\in \mathcal{P}(E_{i,l})$ that make the equality $$\varphi(x) =
\sum_{a\in \mathcal{P}(E_{i,l})} c_a \cdot a_x$$ true for all elements
$x$ as above. We want to prove that the only non-zero terms in this
sum are the ones which correspond to ``concatenation'' of the symbols
$s_k=\{x_{1,k}, \ldots, x_{i,k}\}$.  Precisely, let $A_{i,l}$ be the
subset of $\mathcal{P}(E_{i,l})$ consisting of elements $a$ such that
if $(m,k) \in a$ for some $m$ and $k$, then $(m',k) \in a$ for all
integer $m' \in [1,i]$. Also, let $B_{i,l}$ be the complement of
$A_{i,l}$ in $\mathcal{P}(E_{i,l})$.  Notice that elements in
$A_{i,l}$ have cardinality a multiple of $i$.

\begin{lemma} \label{lemma1} $c_a \neq 0$ implies $a\in A_{i,l}$.
\end{lemma}

\noindent \textbf{Proof of lemma \ref{lemma1}.} We will proceed by
induction on $\# (a)$. Let $P_n$ be the proposition ``For all $a\in
B_{i,l}$ such that $\#(a) \leq n$, $c_a=0$''.

$P_0$ is true since if $\# (a) = 0$, $a$ is not in $B_{i,l}$. Now,
assume $P_n$ is true for some $n$, and let's prove that $P_{n+1}$
holds. Let $a\in B_{i,l}$ be of cardinal $n+1$. Consider the element
$x = \sum_{k=1}^l \{x_{1,k}, \ldots, x_{i,k}\}$ where $x_{m,k}=1$ if
$(m,k) \notin a$. To be precise, this means that we are considering
the image of $\sum_{k=1}^l \{x_{1,k}, \ldots, x_{i,k}\}$ under the
successive application of the specialization maps $s_{x_{m,k}}$ for
$(m,k) \notin a$. Then, by induction hypothesis,

$$\varphi(x) = \sum_{a' \subseteq a} c_{a'} \cdot a'_x = c_a \cdot a_x +
\sum_{\fontsize{7pt}{5}\selectfont{\begin{array}{c} a' \varsubsetneq a \\
      a'\in A_{i,l} \end{array}}} c_{a'}\cdot a'_x.$$ Also, because
$A_{i,l}$ is stable under union, we have
$$\bigcup_{\fontsize{7pt}{5}\selectfont{\begin{array}{c} a'
      \varsubsetneq a \\ a'\in A_{i,l} \end{array}}} a' \varsubsetneq
a.$$ So now, consider an element $(m,k) \in a- \bigg( \bigcup_{a'
  \subset a, a\in A_{i,l}} a' \bigg)$ for some $m$ and $k$, then there
exists an $m'$ such that $(m',k) \notin a$. As $x_{m',k} =1$, we see
that $\varphi(x)$ doesn't depend on $x_{m,k}$. The only term in
$\varphi(x)$ where $x_{m,k}$ appears is $c_a \cdot a_x$. Therefore,
for any field extension $k/k_0$, and for any values taken in
$\KM_1(k)/p$ assigned to the elements $x_{m,k}$ for $(m,k)$ belonging
to $a$, we must have $c_a \cdot a_x = 0$. Proposition \ref{inj}
implies $c_a=0$. \qed \medskip

Hence we obtain $$\varphi(x)=\sum_{a\in A_{i,l}} c_a \cdot a_x.$$

\begin{lemma} \label{lemma2} If $a$, $a' \in A_{i,l}$ are such that
  $\# (a) =\# (a') =ri \in \Z i$, then $c_a =c_{a'} \in
  \KM_{j-ri}(k_0)/p$.
\end{lemma}

\noindent \textbf{Proof of lemma \ref{lemma2}.} It suffices to impose
$l-r$ of the symbols appearing in the decomposition of $x$ to be $0$
and to use the commutativity of addition. \qed\medskip

All in all, noticing that $\sum_{a\in A_{i,l}, \#(a) = ri} a_x =
\gamma_r(x)$, we obtain the existence of elements $y_r \in
\KM_{j-ri}(k_0)/p$ such that for all $x \in \KM_i(k)/p$, $$\varphi(x)
= \sum_{r\geq 0} y_r \cdot \gamma_r(x).$$ \qed \medskip

\begin{proposition}\label{odd} In the case where $i$ and $p$ are odd, the
  algebra of operations is the free $\KM_*(k_0)/p$-module of rank $2$
  generated by $\gamma_0$ and $\gamma_1$. \end{proposition}

\noindent \textbf{Proof.} The module is clearly free and its elements
do define operations. So now, let $\varphi$ be an operation
$\KM_i(k)/p \r \KM_j(k)/p$ over $k_0$ and let $x=s_1 + \ldots + s_l
\in \KM_i(k)/p$ be an element of length at most $l$. With the same
notations as in the proof of Proposition \ref{main1}, we have
$\varphi(x) = \sum_{a\in E_{i,l}} c_a \cdot a_x$ and Lemma
\ref{lemma1} applies. So, actually, $\varphi(x) = \sum_{a\in A_{i,l}}
c_a \cdot a_x $. We want to show that $c_a=0$ as soon as $\#(a) > i$.
It is possible to write
$$\varphi(s_1 + \ldots + s_l) = s_1 \cdot s_2 \cdot \varphi_0(s_3 +
\ldots + s_l) + s_1\cdot \varphi_1(s_3 + \ldots + s_l) + s_2\cdot
\varphi_2(s_3 + \ldots + s_l).$$ If $s_1$ and $s_2$ are permuted, we
should obtain the same result. Substracting both identities and
considering that $s_1 \cdot s_2 = - s_2 \cdot s_1$, we get the
equality
$$2s_1 \cdot s_2 \cdot \varphi_0(s_3 + \ldots + s_l) + (s_1-s_2) \cdot
\big(\varphi_1(s_3 + \ldots + s_l) - \varphi_2(s_3 + \ldots +
s_l)\big) = 0.$$ Setting $s_2=0$ gives $s_1\cdot \big(\varphi_1(s_3 +
\ldots + s_l) - \varphi_2(s_3 + \ldots + s_l)\big) = 0$ for all
$s_1,s_3,\ldots, s_l$ and Proposition \ref{inj} implies $\varphi_1(s_3
+ \ldots + s_l)=\varphi_2(s_3 + \ldots + s_l)$. Hence, we are led to
the equality $$2s_1 \cdot s_2 \cdot \varphi_0(s_3 + \ldots + s_l) =
0$$ for all symbols $s_1,\ldots,s_l$. Since $2$ is invertible in
$\KM_0(k_0)/p = \mathbf{F}_p$, Proposition \ref{inj} implies
$\varphi_0(s_3 + \ldots + s_l) = 0$ for all $s_3,\ldots,s_l$ that is
$\varphi_0=0$, since $l$ was arbitrary. The result follows by taking
all the different pairs of symbols in place of $s_1$ and $s_2$. \qed
\medskip

\subsection{Operations in Milnor $K$-theory mod $2$}

We finish proving Theorem \ref{main2} by considering the remaining
case, that is $p=2$ and $-1$ is not necessarily a square in the base
field. Let $k_0$ be any field and consider again, as in section 2.4,
the map $\tau_i : \KM_*(k_0)/2 \r \KM_*(k_0)/2, \ x \mapsto
\{-1\}^{i-1} \cdot x$.

\begin{proposition}\label{opmod2}  The algebra of operations $\KM_i(k)/2 \r
  \KM_*(k)/2$ over $k_0$ commuting with field extensions is
  \begin{itemize} 
  \item If $i=0$, the free $\KM_*(k_0)/2$-module of rank $2$ of
    functions $\mathbf{F}_2 \r \KM_*(k_0)/2$. 
  \item If $i=1$, the free $\KM_*(k_0)/2$-module of rank $2$,
    generated by $\gamma_0$ and $\gamma_1$.
  \item If $i \geq 2$, the $\KM_*(k_0)/2$-module $$\KM_*(k_0)/2\cdot
    \gamma_0 \oplus \KM_*(k_0)/2\cdot \gamma_1 \oplus \bigoplus_{n\geq
      2} \mathrm{Ker}(\tau_i) \cdot \gamma_n.$$ 
  \end{itemize}
\end{proposition}

\noindent \textbf{Proof.} The cases $i=0$ and $i=1$ have already been
treated. Now, suppose $i \geq 2$. Given $l\geq 2$ and $y_l \in
\mathrm{Ker}(\tau_i)$, the map $y_l \cdot \gamma_l$ is a well-defined
operation by Proposition \ref{mod2}. The inclusion ``$\supseteq$''
holds for the very same reason as in the proof of Proposition
\ref{main1}.

The proof of Proposition \ref{main1} shows that if $\varphi$ is an
operation $\KM_i(k)/2 \r \KM_j(k)/2$, then necessarily there exist
elements $y_r \in \KM_{j-ri}(k_0)/2$ such that $\varphi = y_0 +
y_1\cdot \gamma_1 + \ldots + y_l\cdot \gamma_l + \ldots$. All we have
to prove is that necessarily, for a given integer $l \geq 2$, $y_l$
must satisfy $y_l\cdot \{-1\}^{i-1}=0$. Let $x= s_1 + \ldots + s_l$ be
a sum of $l \geq 2$ symbols. Suppose $s_1 =
\{x_1,\ldots,x_{i-1},x_i\}$ and $s_2=\{x_1,\ldots,x_{i-1},y_i\}$. Then
$$\varphi(x) = y_0 + y_1\cdot x + \ldots + y_l\cdot \gamma_l (s_1 +
s_2 + \ldots + s_l)$$ but also
$$\varphi(x) = y_0 + y_1\cdot x + \ldots + y_{l-1}\cdot \gamma_{l-1} (
\{x_1,\ldots,x_{i-1},x_iy_i\}+ s_3 + \ldots + s_l)$$ It is easy to
check that the difference of these two equalities lead to the equality
$$y_l \cdot \{x_1,\ldots,x_{i-1},x_i,x_1,\ldots,x_{i-1},y_i\}\cdot s_3
\cdot \ldots \cdot s_l = 0.$$ Recall that $\{x,x\}=\{x,-1\}$ for all
$x \in k^\times$ and that $\KM_*(k)/2$ is a commutative algebra. Hence
$y_l$ must satisfy $y_l \cdot \{-1\}^{i-1}\cdot
\{x_1,\ldots,x_{i-1},x_i,y_i\}\cdot s_3 \cdot \ldots \cdot s_l=0$.
Proposition \ref{inj} then implies that necessarily $y_l \cdot
\{-1\}^{i-1} =0$. \qed \medskip

As a consequence of Theorem \ref{main} and Theorem \ref{main2}, we get

\begin{corollary}
  Let $k_0$ be any field and let $i$ be an integer. Then, the
  operations $\varphi : \KM_i/p \longrightarrow \KM_*/p$ over $k_0$
  commute with specialization maps.
\end{corollary} \qed

\subsection{Operations in integral Milnor $K$-theory}

Operations $\KM_i \r \KM_*$ over $k_0$ are not as nice as in the mod
$p$ case since such operations are not determined by the image of a
transcendental element. For example, let $a$ and $b$ be distinct
elements in $\KM_1(k_0)$ (for $k_0$ not the field with only $2$
elements). Consider the operation $\varphi : \KM_1 \r \KM_2$ that
assigns to each element $t$ transcendental over $k_0$ the value
$\{t\}\cdot a$ and to each element $e$ algebraic over $k_0$ the value
$\{e\} \cdot b$.  This is a well-defined operation since for any
extension $k$ of $k_0$, and any transcendental element $t \in k$ and
any algebraic element $e\in k$, we have $\{t\} \neq \{e\}$ in
$\KM_1(k)$. Also, such an operation is not of the form described in
Theorem \ref{theorem} because of Proposition \ref{inj}.

Nonetheless, the image $\varphi(\{t\})$ of any transcendental element
$t$ over $k_0$ determines the image $\varphi(\{u\})$ of any other
transcendental element $u$ over $k_0$, via the obvious isomorphism
$k_0(t) \simeq k_0(u)$. Also, it seems natural to impose the
operations $\varphi : \KM_1 \r \KM_*$ to commute with specialization
maps, in which case the image of any algebraic element over $k_0$ is
determined by $\varphi(\{t\})$.

\begin{definition}
  Let $k_0$ be any field. We say that an operation $\varphi : \KM_i
  \longrightarrow \KM_*$ over $k_0$ commutes with specialization maps
  if $\varphi$ satisfies the conclusion of Definition \ref{spedef}. In
  particular, $\varphi$ commutes with specialization maps only if for
  any extension $k/k_0$, any $t$ transcendental over $k$ and for any
  closed point $P$ in $\P^1_k$, we have a commutative diagram
  \begin{center} $ \xymatrix{\KM_i(k(t)) \ar[r]^{
        \ \varphi} \ar[d]_{s_\pi} & \KM_*(k(t)) \ar[d]_{s_\pi} \\
      \KM_i(\kappa_P) \ar[r]^{ \ \varphi} & \KM_*(\kappa_P) }$
  \end{center} where $\kappa_P$ denotes the residue field of $k(t)$
  with respect to the valuation $v_P$ corresponding to the polynomial
  $P$, and $\pi$ is any uniformizer for the valuation $v_P$.
\end{definition}

Before we describe operations commuting with specialization maps, we
need two lemmas. Firstly, residue maps and specialization maps are
well-behaved with respect to transcendental field extensions.

\begin{lemma} \label{com} For any field $F$ and any $u$ transcendental
  over $F$, let $\iota_u : \KM_*(F) \r \KM_*(F(u))$ be the injective
  map induced by the inclusion of field $F \subset F(u)$. Let $k$ be a
  field and $P$ a closed point in the projective line $\P^1_k$ with
  residue field $\kappa_P$. Let $v_P$ be the valuation on $k(t)$
  corresponding to $P$ and let $\pi$ be a local parameter for $v_P$ in
  $k(t)$.

  Then, the valuation $v_P$ extends naturally to a valuation, that we
  still write $v_P$, on $k(u)(t)$, and $\pi$ seen as an element in
  $k(u)(t)$ defines a uniformizer for $v_P$ in $k(u)(t)$. Moreover,
  the residue map $\partial_{v_P}$ and the specialization map $s_\pi$
  commute with $\iota_u$. Precisely, the following diagrams commute :
  \begin{center} $ \xymatrix{\KM_*(k(t)) \ar[d]_{\partial_{v_P}}
      \ar[r]^{\iota_u \ } & \KM_*(k(u,t)) \ar[d]_{\partial_{v_P}} \\
      \KM_{*-1}(\kappa_P) \ar[r]^{\iota_u \ } &
      \KM_{*-1}(\kappa_P(u))}$ $\ \ \ \ \ \ \ \ $ $
    \xymatrix{\KM_*(k(t)) \ar[d]_{s_\pi} \ar[r]^{\iota_u \ } &
      \KM_*(k(u,t)) \ar[d]_{s_\pi} \\ \KM_{*}(\kappa_P)
      \ar[r]^{\iota_u \ } & \KM_{*}(\kappa_P(u)) }$

   \end{center}
  
\end{lemma}

\noindent \textbf{Proof.} The fact that $v_P$ and $\pi$ extend
respectively to a valuation and to a uniformizer on $k(u)(t)$ is
straightforward. The commutativity of the diagrams is an immediate
consequence of the definition of the residue map and of the
specialization map. \qed \medskip

Secondly, We need to relate the specialization maps $s_\pi$ and
$s_{\pi'}$ for two different choices of uniformizers $\pi$ and $\pi'$.
If $P$ is a closed point in $\P^1_{k}-\{\infty\}$ and $v_P$ is the
corresponding valuation on $k(t)$, then $P$ is a local parameter for
$v_P$. Now, if $Q \in k(t)$ is such that $v_P(Q)=0$, then $PQ$ defines
another local parameter for $v_P$. Therefore, it is possible to
consider specialization maps $s_P$ and $s_{PQ}$ mapping $\KM_*(k(t))$
to $\KM_{*-1}(\kappa_P)$.

\begin{lemma} \label{spe} Let $k$ be a field and $x$ be an element in
  $\KM_*(k(t))$. If $P$ is a closed point in $\P^1_{k}-\{\infty\}$ and
  $Q \in k(t)$ is such that $v_P(Q)=0$, then we have the formula
  $$s_{PQ}(x) = s_P(x) - s_P\big(\{-Q(t)\}\big) \cdot
  \partial_{v_P}(x).$$
\end{lemma}

\noindent \textbf{Proof.} Under the assumption made on $Q$, the
element $PQ$ is a uniformizer for the valuation $v_P$. Hence, we have
$s_{PQ}(x)=\partial_{v_P}\big( \{-PQ\}\cdot x \big) = s_P(x) +
\partial_{v_P}\big( \{-Q\}\cdot x \big)$. It is thus enough to show
that $\partial_{v_P} \big( \{-Q\} \cdot x \big) = -s_P \big(
\{-Q\}\big) \cdot \partial_{v_P}(x)$. The element $Q$ being a unit in
the ring $\{a \in k(t), \ v_P(a) \geq 0\}$, this follows from the very
definition of the residue and specialization maps. \qed \medskip

\begin{proposition} Let $k_0$ be a field. \begin{itemize}
  \item The algebra of operations $\KM_0(k) \r \KM_*(k)$ over $k_0$
    commuting with specialization maps is the $\KM_*(k_0)$-module
    of functions $\Z \r \KM_*(k_0)$.
  \item The algebra of operations $\KM_1(k) \r \KM_*(k)$ over $k_0$
    commuting with specialization maps is the free $\KM_*(k_0)$-module
    generated by $\gamma_0$ and $\gamma_1$.
\end{itemize} 
\end{proposition}

\noindent \textbf{Proof.} For the first statement, $\KM_0$
is the constant functor with value $\Z$. Also, inclusion of fields and
specialization maps induce the identity on the $\KM_0$-groups of
fields. Hence the result.

For the second statement, such an operation $\varphi : \KM_1(k) \r
\KM_*(k)$ is determined by the image of $\{t\} \in \KM_1(k_0(t))$ as
discussed above. Let $\varphi : \KM_1 \r \KM_*$ be an operation over
$k_0$ commuting with specialization maps. We are going to show that,
for $t$ a transcendental element over $k_0$, $\varphi(\{t\})$ is
unramified outside $\{0,\infty\}$. By Milnor's exact sequence, this
will prove the Proposition.

Let's consider the function field in one indeterminate $k=k_0(u)$ and
a monic irreducible polynomial $P \in k_0[t]$ that we can also see as
a monic irreducible polynomial in $k[t]$ with coefficients in $k_0$,
via the obvious inclusion of $k_0$ into $k$. Let $v_P$ denote the
valuation on $k(t)$ corresponding to $P$. The choice of a uniformizer
$\pi \in k(t)$ is equivalent to the choice of an element $Q \in k(t)$
such that $v_P(Q)=0$, by setting $\pi_Q=PQ$. Let $\alpha$ be the image
of $t$ in the residue field $\kappa_P = k[t]/P$. By definition of
$\varphi$, we have a commutative diagram for any $Q \in k_0(u,t)$ such
that $v_P(Q)=0$
\begin{center} $\xymatrix{ \KM_1(k_0(u,t)) \ar[r]^{ \varphi}
    \ar[d]_{s_{\pi_Q}}
    &   \KM_*(k_0(u,t))  \ar[d]_{s_{\pi_Q}} \\
    \KM_1(\kappa_P(u)) \ar[r]^{ \varphi } & \KM_*(\kappa_P(u)) }$
\end{center}
For any field $F$, let $\iota_u : \KM_*(F) \r \KM_*(F(u))$ be the
injective map induced by the inclusion of field $F \subset F(u)$. For
$\{t\} \in \KM_1(k_0(t))$, we then have $s_{\pi_Q}\circ \varphi \circ
\iota_u (\{t\}) = \varphi \circ s_{\pi_Q} \circ \iota_u (\{t\})$. If
$P \neq t$, Lemma \ref{spe} says that, on the one hand, $\varphi \circ
s_{\pi_Q} \circ \iota_u (\{t\}) = \iota_u
\big(\varphi(\{\alpha\})\big)$ and on the other hand $s_{\pi_Q}\circ
\varphi \circ \iota_u (\{t\}) = s_{\pi_Q}\circ \iota_u \circ \varphi
(\{t\}) = \{-Q^{-1}(\alpha)\} \cdot
\partial_{v_P} \big( \iota_u \circ \varphi (\{t\})\big) + s_P\big(
\iota_u(\varphi(\{t\})) \big)$. Also, since $P$ has its coefficients in
$k_0$, Lemma \ref{com} implies that $\partial_{v_P} \big( \iota_u
\circ \varphi (\{t\})\big) = \iota_u \big(\partial_{v_P} \circ \varphi
(\{t\})\big)$ and also that $ s_P \big(
\iota_u \circ \varphi(\{t\}) \big)= \iota_u\big(
s_P \circ \varphi(\{t\}) \big)$. All in all, we have
$$\iota_u \big(\varphi(\{\alpha\}) \big) =  \{-Q^{-1}(\alpha)\} \cdot
\iota_u \big(\partial_{v_P} \circ \varphi (\{t\})\big) + \iota_u\big(
s_P(\varphi(\{t\})) \big) \in \KM_*(\kappa_P(u)).$$ Basically, we have
just been expressing the fact that apart from $Q$, everything exists
before adjoining that indeterminate $u$. So now, let $Q$ be the
constant polynomial in $k_0(u)[t]$ equal to $-u^{-1}$ and consider the
residue map $\partial_u : \KM_*(\kappa_P(u)) \r \KM_{*-1}(\kappa_P)$.
Applying $\partial_u$ to the above equality, and using the injectivity
of $\iota_u$, we get
$$\partial_{v_P} \circ \varphi (\{t\})=0.$$ \qed \medskip

\begin{proposition} \label{integral} Let $k_0$ be any field and
$i$ an integer $\geq 2$. Write $\tau_i : \KM_*(k_0) \r \KM_*(k_0)$ for the homomorphism $ x \mapsto \{-1\}^{i-1} \cdot x$.
\begin{itemize} 
\item If $i$ is even, the algebra of operations $\KM_i(k) \r \KM_*(k)$
  over $k_0$ commuting with specialization maps is the
  $\KM_*(k_0)$-module $$\KM_*(k_0) \oplus \KM_*(k_0)\cdot\mathrm{id}
  \oplus \bigoplus_{n\geq 2} \mathrm{Ker}(\tau_i)\cdot \gamma_n.$$
\item If $i$ is odd, the algebra of operations $\KM_i(k) \r \KM_*(k)$
  over $k_0$ commuting with specialization maps is the
  $\KM_*(k_0)$-module $$\KM_*(k_0) \oplus \KM_*(k_0)\cdot\mathrm{id}
  \oplus \bigoplus_{n\geq 2} {}_2\mathrm{Ker}(\tau_i)\cdot \gamma_n.$$
 \end{itemize} 
\end{proposition}

\noindent \textbf{Proof.} The proof is the same as for Theorems
\ref{main} and \ref{main2}. An operation $\varphi : \KM_i \r \KM_*$
over $k_0$ is necessarily a sum of divided power operations and it is
well-defined if and only if these are weak divided power operations as
in section 2.5. \qed \medskip

When $\mathrm{char} \ k_0 = 2$, the maps $\tau_i$ are zero for $i \geq
2$. Also, in \cite{Izhboldin}, Izhboldin proves that if $k_0$ has
characteristic $p$ then the Milnor $K$-groups $\KM_n(k_0)$ have no
$p$-torsion (result conjectured by Tate). Hence, when $\mathrm{char} \
k_0 = 2$, ${}_2 \KM_*(k_0) = 0$, the above Proposition becomes

\begin{proposition} Let $k_0$ be a field of characteristic $2$ and $i$ an
integer $\geq 2$. 
\begin{itemize}
\item If $i$ is even, the algebra of operations $\KM_i(k) \r \KM_*(k)$
  over $k_0$ commuting with specialization maps is the free
  $\KM_*(k_0)$-module $$\KM_*(k_0) \oplus \KM_*(k_0)\cdot\mathrm{id}
  \oplus \bigoplus_{n\geq 2} \KM_*(k_0)\cdot \gamma_n.$$
\item If $i$ is odd, the algebra of operations $\KM_i(k) \r \KM_*(k)$
  over $k_0$ commuting with specialization maps is the free
  $\KM_*(k_0)$-module $$\KM_*(k_0) \oplus
  \KM_*(k_0)\cdot\mathrm{id}.$$
\end{itemize}

\end{proposition} \qed

\subsection{additive operations}

An additive operation $\varphi : \KM_i \r \KM_*$ over $k_0$ is an
operation which enjoys the following property : for all field $k/k_0$
and for all $x$ and $y$ in $\KM_i(k)$, $\varphi(x+y)=\varphi(x) +
\varphi(y)$. The set of additive operations over $k_0$ has clearly the
structure of a $\KM_*(k_0)$-module.

For example, an operation in mod $p$ Milnor $K$-theory is a sum of
divided power operations, and from there it is possible to prove that
an additive operation $\KM_i/p \r \KM_*/p$ is necessarily of the form
$x \mapsto a \cdot x$ for some $a \in \KM_*(k_0)/p$. As already
mentioned in subsection 3.4, operations $\KM_i \r \KM_*$ over $k_0$
can be pathological. However, additive operations in integral Milnor
$K$-theory have a nice description (and the mod $p$ case can be proven
the same way) :

\begin{proposition} \label{additive} Let $k_0$ be a field. The algebra
  of additive operations $\varphi : \KM_i \r \KM_*$ over $k_0$ is the
  free $\KM_*(k_0)$-module generated by the identity map. In other
  words, given an additive operation $\varphi$ over $k_0$, there
  exists $a \in \KM_*(k_0)$ such that for all field $k/k_0$ and all
  $x\in \KM_i(k)$, $\varphi(x)=a\cdot x$.
\end{proposition}

\noindent \textbf{Proof.} We start by proving the case $i=1$. Let
$\varphi : \KM_1 \r \KM_*$ be an additive operation over $k_0$. The
proof is very much the same as the proof of Theorem 3.4. We first
claim that $\varphi$ is determined by $\varphi(\{t\})$ for $t$ a
transcendental element over $k_0$. Indeed, if $u$ is another
transcendental element over $k_0$, the isomorphism $k_0(t) \simeq
k_0(u)$ mapping $t$ to $u$ determines $\varphi(\{u\})$. If $e$ is an
algebraic element over $k_0$, then $et$ is transcendental over $k_0$
and $\varphi(\{e\}) = \varphi(\{et\}) - \varphi(\{t\})$ in
$\KM_*(k_0(e,t))$. By Milnor's exact sequence, $\KM_*(k_0(e))$ embeds
into $\KM_*(k_0(e,t))$ and thus $\varphi(\{e\})$ is uniquely
determined as an element of $\KM_*(k_0(e))$.

Therefore, it is enough to show that $\varphi(\{t\}) \in
\KM_*(k_0(t))$ has possibly non-zero residue only at the points $0$
and $\infty \in \P_{k_0}^1$. For this purpose, let $t$ and $u$ be two
algebraically independent transcendental elements over $k_0$. By
additivity, we have $\varphi(\{ut\}) = \varphi(\{u\}) +
\varphi(\{t\})$ in $\KM_*(k_0(u,t))$. Let $P$ be a point in
$\P_{k_0}^1$, i.e. a monic irreducible polynomial with coefficients in
$k_0$. Let's write $P_u$ for the polynomial $P_u(t)=P(ut) \in
k_0(u)[t]$. The same arguments as in Step 2 of the proof of theorem
3.4 show that $P_u$ must be equal to $\alpha P$ for some $\alpha \in
k_0(u)$. This implies that $P$ must be equal to $t$. Moreover, if $c$
is the specialization at infinity of $\varphi(\{t\})$, the formula $c
= s_\infty \varphi(\{ut\}) = s_\infty \varphi(\{u\}) + s_\infty
\varphi(\{t\}) = 2c$ shows that $\varphi$ has vanishing specialization
at infinity. Therefore, if $a= \partial_0 \varphi(\{t\}) \in
\KM_*(k_0)$, we have for all field extension $k/k_0$ and all $x \in
\KM_1(k)$, $\varphi(x)= a \cdot x$.  This clearly defines an additive
operation.

We now finish the proof by induction on $i$. Assume all additive
operations $\KM_{i-1} \r \KM_*$ over $k_0$ are of the form $a\cdot
\mathrm{id}$ for some $a \in \KM_*(k_0)$ and consider an operation
$\varphi : \KM_i \r \KM_*$. By additivity, it is enough to determine
the image of symbols. Let $k$ be a field extension of $k_0$ and
$x_1,\ldots, x_{i-1}$ be elements in $k^\times$. Let $F$ be a field
extension of $k$ and $x \in F^\times$. The map $x \mapsto
\varphi(\{x_1,\ldots,x_{i-1},x\})$ defines an additive operation $
\KM_1 \r \KM_*$ over $k$. Hence, there exists $a_{x_1,\ldots, x_{i-1}}
\in \KM_*(k)$ such that $\varphi(\{x_1,\ldots,x_{i-1},x\}) =
a_{x_1,\ldots, x_{i-1}} \cdot \{x\}$ for all extension $F/k$ and all
$x \in F^\times$. Now, it is easy to check (thanks to Prop. 3.1 for
example) that the map $\{x_1,\ldots,x_{i-1}\} \mapsto
a_{x_1,\ldots,x_{i-1}}$ induces an additive operation $ \KM_{i-1} \r
\KM_*$ over $k_0$. We conclude with the induction hypothesis. \qed

\section{Operations in Milnor $K$-theory of a smooth scheme}

In this section, we generalize the results about operations in Milnor
$K$-theory of fields to the case of smooth schemes over a field $k$.
We are first interested in the Milnor $K$-theory of a regular
$k$-scheme defined as the kernel of the first map in the Gersten
complex. Such a definition coincides with Rost's Chow groups with
coefficients as constructed in \cite{Rost}. Indeed, for $X$ a regular
$k$-scheme of dimension $d$, with Rost's notations, we have $\KM_n(X)
= A_d(X,n-d)$ where $A_p(X,q)$ is the $p^\mathrm{th}$ homology group
of the Gersten complex $C_*(X,q)$ defined by $C_p(X,q) = \bigoplus_{x
  \in X_{(p)}} \KM_{p+q}\big(k(x)\big)$ and $X_{(p)}$ denotes the
$p$-dimensional points in $X$. It is then a fact that $\KM_n$ defines
a contravariant functor from the category of smooth $k$-schemes to the
category of groups. As in the case of fields, we are able to determine
all operations $\KM_n/p \r \KM_*/p$ over a smooth $k$-scheme $X$. In
view of the Gersten complex, we can write $\KM_n(X) = H^0(X, \Km_n)$.
It is then possible, under the assumption of the Bloch-Kato
conjecture, to relate for $p \neq \mathrm{char} \ k$ the Milnor
$K$-group $\KM_n(X)/p$ and the unramified cohomology group
$H^0\big(X,\mathcal{H}^n_{et}(\Z/p)\big)$, and thus to describe all
the operations on the unramified cohomology of smooth schemes over
$k$.

We are then interested in the Milnor $K$-theory
$\bar{K}^{\mathrm{M}}_n(A)$ of a ring $A$ defined as the tensor
algebra of the units in $A$ subject to the Steinberg relations.  This
defines a covariant functor from the category of rings to the category
of sets. If $k$ is an infinite field and if $A$ is a regular
semi-local $k$-algebra, we are also able to determine all operations
$\bar{K}^{\mathrm{M}}_n/p \r \bar{K}^{\mathrm{M}}_*/p$ over $A$.

\subsection{The unramified case} 

Let $X$ be a regular (in codimension-$1$) scheme, and denote by
$X^{(r)}$ the set of codimension-$r$ points in $X$. If $x$ is a
codimension-$0$ point of $X$, e.g. the generic point of $X$ if $X$ is
irreducible, then the codimension-$1$ points in the closure of $x$
define discrete valuations on the function field $k(x)$ of $x$, and
thus residue maps on the Milnor $K$-theory of $k(x)$. We define the
Milnor $K$-theory of the scheme $X$ to be
  $$\KM_n(X) = \ker \Big( \bigoplus_{x \in X^{(0)}}
  \KM_n\big(k(x)\big) \stackrel{\partial}{\longrightarrow}
  \bigoplus_{y \in X^{(1)}} \KM_{n-1}\big(k(y)\big) \Big).$$ In
  particular, this definition makes sense for regular rings. For a
  regular scheme $X$ assumed to be irreducible and for any $i\geq 1$,
  an element $x \in \KM_i(X)$ is an element of $\KM_i\big(k(X)\big)$
  which is unramified along all codimension-$1$ points of $X$, i.e.
  which has zero residue for all residue maps corresponding to
  codimension-$1$ points in $X$. We say that an element $x \in
  \KM_i\big(k(X)\big)$ is unramified if it belongs to $\KM_i(X)$. It
  is therefore possible to write an element $x$ of $\KM_i(X)$ as a sum
  of symbols $s_k = \{x_{1,k},\ldots,x_{i,k}\}$, $1 \leq k \leq l$,
  with all the $x_{j,k}$'s in $k(X)^\times$. Given an integer $n$ and
  an element $y_n \in \KM_*(X)$, an $n^\mathrm{th}$ divided power of
  $x$ written as a sum $\sum_{k=1}^l s_k$ is
$$y_n \cdot \gamma_n(x) = y_n \cdot \sum_{1
  \leq l_1 < \ldots < l_n \leq l} s_{l_1} \cdot \ldots \cdot
s_{l_n}.$$

\begin{lemma}\cite[Prop 7.1.7]{Gille} \label{unram}
  Let $K$ be a field equipped with a discrete valuation $v$, and let
  $\mathcal{O}_K$ be its ring of integers and $\kappa$ be its residue
  field. Then, $\ker \big(\KM_n(K)
  \stackrel{\partial_v}{\longrightarrow} \KM_{n-1}(\kappa) \big)$ is
  generated as a group by symbols of the form $\{x_1,\ldots,x_n\}$
  where the $x_i$'s are units in $\mathcal{O}_K$ for all $i$.
\end{lemma}

Let's mention that this lemma implies the following.

\begin{corollary}
  Let $X$ be a regular scheme. The cup-product on
  $\KM_*\big(k(X)\big)$ endows the group $\KM_*(X) = \bigoplus_{n \geq
    0} \KM_n(X)$ with a ring structure.
\end{corollary}

\begin{proposition} \label{dprings} If $X$ is a regular scheme, divided powers
  are well-defined on $\KM_i(X)/p$ in the following cases :
  \begin{itemize}
  \item if $i=0$ or $i=1$.
  \item if $p=2$, $i\geq 2$ and $y_n \in \ker \big( \tau_i : \KM_*(X)/2 \r
    \KM_*(X)/2, \ x \mapsto \{-1\}^{i-1} \cdot x \big)$.
  \item if $p$ is odd and $i\geq 2$ is even, and $y_n$ is any element
    in $\KM_*(X)/p$.
  \end{itemize}

\end{proposition}

\noindent \textbf{Proof.} For simplicity, assume $X$ is irreducible
with field of rational functions $k(X)$. By the results of section 2,
it suffices to check that if $x = \sum_{k=1}^l s_k \in
\KM_i\big(k(X)\big)$ is unramified then $\gamma_n(x) := \sum_{1 \leq
  l_1 < \ldots < l_n \leq l} s_{l_1} \cdot \ldots \cdot s_{l_n} \in
\KM_i\big(k(X)\big)$ is unramified. So, let $y$ be a codimension-$1$
point in $X$ with local ring $\mathcal{O}_{X,y}$, and let $a$ and $b$
be symbols in $\KM_i\big(k(X)\big)$ unramified along $y$. Then, thanks
to Lemma \ref{unram}, we can write $a$ as a sum of symbols
$\{a_1,\ldots,a_n\}$ and $b$ as a sum of symbols $\{b_1,\ldots,b_n\}$
with $a_1,\ldots,a_n,b_1,\ldots,b_n$ being units in
$\mathcal{O}_{X,y}$. It is then clear that $a\cdot b$ is also
unramified along $y$ which finishes the proof. \qed \medskip

The definition of Milnor $K$-theory we gave is functorial with respect
to open immersions of regular schemes. Indeed, if $U \hookrightarrow
X$ is an open immersion of regular scheme, the group homomorphism
$\KM_n(X) \r \KM_n(U)$ is just defined by restriction. Indeed,
divisors on $U$ map injectively into the set of divisors on $X$, and
thus an element $x$ unramified along divisors in $X$ will surely be
unramified along divisors in $U$. We define the Zariski sheaf $\Km_*$
on $X$ to be $$U \mapsto \KM_*\big(U\big)$$ for any Zariski open
subset $U$ of $X$.  Clearly, we have $\KM_i(X) = H^0(X,\Km_i)$.  By a
map of sheaves $\varphi : \Km_i/p \r \Km_*/p$, we mean a map that
commutes with open immersions, i.e. if $U \hookrightarrow X$ is an
open immersion of regular schemes, the following diagram commutes
\begin{center} $ \xymatrix{\KM_i(X) \ar[r]^{
      \ \varphi} \ar[d] & \KM_*(X) \ar[d] \\
    \KM_i(U) \ar[r]^{ \ \varphi} & \KM_*(U) }$
  \end{center}As a straightforward consequence of the
above, we have

\begin{proposition}
  Let $p$ be a prime number and  $X$ be a regular
  scheme.  Then, there exist divided powers of sheaves of sets on
  $X$ $$\gamma_n : \Km_i/p \r \Km_*/p$$ in the following cases :
  \begin{itemize}
    \item if $i=0$ or $i=1$.
    \item if $p=2$, $i\geq 2$ and $-1$ is a square in $\mathcal{O}_X(X)$.
    \item if $p$ is odd and $i\geq 2$ is even.
    \end{itemize}
\end{proposition}

\noindent \textbf{Proof.} This is clear from the definitions and
Proposition \ref{dprings}. \qed \medskip

\begin{remark}
  In the case when $p=2$, $i\geq 2$ and $-1$ is not a square in
  $\mathcal{O}_X(X)$, it is still possible to define some operations
  $\Km_i/2 \r \Km_*/2$.  Indeed, if $\tau_i :
  \KM_*(\mathcal{O}_X(X))/2 \r \KM_*(\mathcal{O}_X(X))/2, \ x \mapsto
  \{-1\}^{i-1}\cdot x$ and if $y_n \in \ker \tau_i$, then we have an
  operation of sheaves on $X$, $y_n \cdot \gamma_n : \Km_i/2 \r
  \Km_*/2$.
\end{remark}

We have seen that $\KM_n$ is functorial with respect to open
immersions of regular schemes. If $k$ is a field, it is actually
functorial with respect to any map between smooth $k$-schemes. Given a
map $f : Y \r X$ between smooth $k$-schemes, Rost constructs in
\cite[section 12]{Rost} a pull-back group homomorphism $f^* : \KM_n(X)
\r \KM_n(Y)$, and shows that it is functorial. In particular, if $f :
Y \r X$ is a dominant map of smooth $k$-schemes, it induces an
embedding of the field of functions $k(X)$ of $X$ into the function
field $k(Y)$ of $Y$, and the map $f^* : \KM_n(X) \r \KM_n(Y)$ is
induced by the map $i : \KM_n\big(k(X)\big) \r \KM_n\big(k(Y)\big)$
coming from the inclusion of fields $k(X) \hookrightarrow k(Y)$
(\cite[Lemma 12.8.]{Rost}), so that an unramified element of
$\KM_n\big(k(X)\big)$ will map to an unramified element of
$\KM_n\big(k(Y)\big)$ under $i$.

Let $k$ be a field and $X$ be a smooth $k$-scheme. Let's denote by
$\mathsf{Sm}_X$ the category of smooth $k$-schemes with a morphism to
$X$ and with morphisms being morphisms of $k$-schemes respecting the
$X$-structure, i.e.  commutative
diagrams  \begin{center} $\xymatrix{ Y \ar[rr] \ar[dr] && Z \ar[dl] \\
    & X }$
\end{center} In particular, if $X$ is irreducible, the spectrum of its
field of rational functions belongs to $\mathsf{Sm}_X$.  The map
$\KM_n/p : \mathsf{Sm}_X \r \mathsf{Sets}$ is a contravariant functor,
and we define an operation over a smooth scheme $X$ to be a natural
transformation from the functor $\KM_n/p : \mathsf{Sm}_X \r
\mathsf{Sets}$ to the functor $\KM_*/p : \mathsf{Sm}_X \r \mathbf{F}_p
- \mathsf{Algebras}$.  Under these assumptions, all the results
concerning fields translate to the case of smooth $k$-schemes and it
is possible to describe all such operations. First, we show that
divided powers commute with Rost's pullback map $f^*: \KM_n(X) \r
\KM_n(Y)$, for a morphism $f:Y \r X$.

\begin{lemma} \label{1} Let $X$ and $Y$ be smooth schemes over a field
  $k$, and let $f : Y \r X$ be a morphism. The pullback map $f^*:
  \KM_n(X)/p \r \KM_n(Y)/p$ commutes with the divided powers of
  Proposition \ref{dprings}.
\end{lemma}

\noindent \textbf{Proof.} The morphism $f : Y \r X$ factors through $Y
\stackrel{i}{\longrightarrow} Y \times X
\stackrel{\pi}{\longrightarrow} X$, where $i$ is the closed immersion
$i(y)=(y,f(y))$ and $\pi$ is the projection $\pi(y,x)=x$. By
functoriality, we have $f^*= i^* \circ \pi^*$. Divided powers commute
with $\pi^*$ because $\pi$ is a dominant map and as such induces a map
on Milnor $K$-theory coming from an inclusion of fields, see above. It
remains to show that divided powers commute with the pullback maps
induced by closed immersions. Let $i : Z \hookrightarrow X$ be a
closed immersion of codimension $c$. By \cite[Corollary 12.4.]{Rost},
$i^* : \KM_n(X) \r \KM_n(Z)$ is the restriction of a composition of
specialization maps $s_1 \circ \ldots \circ s_c : \KM_n\big(k(X)\big)
\r \KM_n\big(k(Z)\big)$. But we already know that divided powers
commute with specialization maps.  \qed \medskip

\begin{theorem} \label{oprings} Let $k$ be any field, $p$ be a prime
  number and $X$ be a smooth scheme over $k$. Operations $\KM_i/p \r
  \KM_*/p$ over the smooth $k$-scheme $X$ are spanned as a
  $\KM_*(X)/p$-module by the divided power operations of Proposition
  \ref{dprings}.
\end{theorem}

\noindent \textbf{Proof.} Lemma \ref{1} shows that divided powers are
indeed operations over $X$ and that so are qny element in their
$\KM_*(X)/p$-span.

For simplicity, assume $X$ is irreducible.  An operation $\varphi :
\KM_i/p \r \KM_*/p$ over $X$ induces naturally an operation
$\bar{\varphi}$ over the field $k(X)$ of rational functions on $X$. By
Theorems \ref{main} and \ref{main2}, the operation $\bar{\varphi}$ is
a sum of divided power operations with coefficients in
$\KM_*\big(k(X)\big)/p$. For any irreducible smooth scheme $Y$ with
field of rational functions $k(Y)$, let's write $\iota$ for the
inclusion of $\KM_*(Y)/p$ inside $\KM_*\big(k(Y)\big)/p$. There exist
elements $y_0,\ldots,y_n$ in $\KM_*\big(k(X)\big)/p$ such that for all
smooth scheme $Y$ over $X$ and for all $x \in \KM_i(Y)/p$, we have
$\varphi \circ \iota (x) = \sum_{k=0}^n y_k \cdot \gamma_k(x)$.  We
also have, by definition of an operation, $\varphi \circ \iota (x)=
\iota \circ \varphi(x)$ for all $x$. If we can prove that the $y_k$'s
are actually in $\KM_*(X)/p$, then we will be done.  First, $y_0$ is
indeed in $\KM_*(X)/p$.  This is because $\varphi(0) = y_0$ must be in
$\KM_*(X)/p$. Suppose we have shown that $y_0,\ldots,y_{l-1}$ are in
$\KM_*(X)/p$ and let's show that $y_l$ is in $\KM_*(X)/p$. Let $Y =
\mathrm{Spec} \ \mathcal{O}_X[t_{j,k}]_{1\leq j \leq i, 1\leq k \leq
  l}$ be the smooth scheme $X \times \mathbf{A}^{il}$. Then, the field
of rational functions of $Y$ is $k(X)(t_{j,k})_{1\leq j \leq i, 1\leq
  k \leq l}$ and if $x = \sum_{k=1}^l \{t_{1,k},\ldots,t_{i,k}\}$,
$\varphi(x) = y_0 + y_1 \cdot x + \ldots + y_l\cdot \gamma_l(x)$.
Therefore, $y_l\cdot \gamma_l(x) = y_l \cdot \{t_{1,1},\ldots,t_{i,1},
\ldots, t_{1,l},\ldots,t_{i,l}\}$ must be in $\KM_*(Y)/p$.

\noindent Also, the closed subschemes $Z$ in $X$ correspond
bijectively to the closed subschemes of the form $Z \times_X Y$ in
$Y$.  Let $u$ be a codimension-$1$ point in $X$ with residue field
$k(u)$, then $u$ corresponds to the codimension-$1$ point $v = u
\times_X Y$ in $Y$.  The residue at $v$ of $y_l\cdot \gamma_l(x)$
considered as an element of $\KM_*\big(k(Y)\big)/p$ is $0$ since
$y_l\cdot \gamma_l(x) \in \KM_*(Y)/p$. We also have
$\partial_v\big(y_l \cdot \gamma_l(x)\big) =
\partial_u(y_l) \cdot \gamma_l(x)$ in $\KM_*\big(k(u)(t_{j,k})_{1\leq
  j \leq i, 1\leq k \leq l}\big)/p$. This in turn implies, by
Proposition \ref{inj}, that $\partial_u(y_l)=0 \in
\KM_*\big(k(u)\big)/p$. Thus, by definition of the Milnor $K$-theory
of a scheme, we get $y_l \in \KM_*(X)/p$.\qed \medskip

\begin{remark}
  Actually, if $X$ is a regular scheme over $k$ and if $\pi : X \times
  \mathbf{A}^r \r X$ is the first projection or more generally if
  $\pi$ is an affine bundle over $X$, then the induced homomorphism
  $\pi^* : \KM_n(X) \r \KM_n(X\times \mathbf{A}^r)$ is an
  isomorphism. See \cite[Prop. 8.6.]{Rost}.
\end{remark}

Let $k$ be a field and $p$ a prime number different from the
characteristic of $k$. We define $\mathcal{H}^i(\Z/p)$ to be the
Zariski sheaf on the category of smooth schemes over $k$ corresponding
to the Zariski presheaf $U \mapsto H^i_{et}\big(U,\Z/p(i)\big)$. If
$X$ is a smooth scheme over $k$, the unramified cohomology of $X$ is
defined to be $H^0\big(X,\mathcal{H}^i(\Z/p)\big)$. This group is
birationally invariant when $X$ is proper over $k$, see Theorem 4.1.1
and Remark 4.1.3 of \cite{CT}. It is worth saying that our unramified
cohomology group is not the same as the one considered by
Colliot-Th\'el\`ene in \cite{CT}. The unramified cohomology is then
clearly functorial with respect to morphisms of smooth $k$-schemes.
Under the assumption of the Bloch-Kato conjecture, the sheaf $\Km_i/p$
maps isomorphically to the sheaf $\mathcal{H}^i(\Z/p)$. Indeed, there
is a morphism of exact Gersten complexes
\begin{center} $\xymatrix{ 0 \ar[r] & \Km_i/p \ar[r] \ar[d] &
    \bigoplus_{x\in X^{(0)}} \KM_i\big(k(x)\big)/p \ar[r]^{\partial}
    \ar[d] & \bigoplus_{x\in X^{(1)}} \KM_{i-1}\big(k(x)\big)/p \ar[d]
    \\ 0 \ar[r] & \mathcal{H}^i(\Z/p) \ar[r] & \bigoplus_{x\in
      X^{(0)}} H^i_{et}\big(k(x), \Z/p(i) \big) \ar[r]^{\partial \ } &
    \bigoplus_{x\in X^{(1)}} H^{i-1}_{et}\big(k(x), \Z/p(\small{i-1}) \big)} $
\end{center}
where the bottom residue map has a description in terms of Galois
cohomology as the edge homomorphisms of a Hochschild-Serre spectral
sequence. It is a fact (see e.g. \cite[section 6.8]{Gille}) that both
residue maps are compatible with the Galois symbol. Hence the claimed
isomorphism of sheaves. In particular, both sheaves have same global
sections, i.e. there is an isomorphism $$\KM_i(X)/p
\stackrel{\simeq}{\longrightarrow}
H^0\big(X,\mathcal{H}^i(\Z/p)\big).$$ Moreover, it can be shown that
this isomorphism is compatible with the pull-back map $f^*$ induced by
any morphism $f:X \r Y$ between smooth $k$-schemes. This proves

\begin{theorem}
  Let $k$ be a field and $p$ be a prime number different from the
  characteristic of $k$. The algebra of operations in unramified
  cohomology $H^0\big(-,\mathcal{H}^i(\Z/p)\big) \r
  H^0\big(-,\mathcal{H}^*(\Z/p)\big)$ is spanned as a
  $H^*_{et}\big(k,\Z/p(*)\big)$-module by the divided powers of
  Proposition \ref{dprings}.
\end{theorem}

\subsection{Operations in Milnor $K$-theory of rings}

A natural question is to ask if $\KM_n(A)$, for a ring $A$, can be
presented with generators and relations. This is optimistic for a
general $A$. However, notice that for a domain $A$, the natural map
$(A^\times)^{\otimes n} \r \KM_n(F)$ factors through
$$\bar{K}^{\mathrm{M}}_n(A) =_{\mathrm{def}} \Big(A^\times \otimes \ldots
\otimes A^\times \Big) / St_A,$$ where $St_A$ is that ideal in
$A^\times \otimes \ldots \otimes A^\times $ generated by elements of
the form $a \otimes (1-a)$ with $a, 1-a \in A^\times$. A ring
homomorphism $A\r B$ induces a ring homomorphism
$\bar{K}^{\mathrm{M}}_*(A) \r \bar{K}^{\mathrm{M}}_*(B)$ and this is
functorial.

If $A$ is an excellent ring, there is a Gersten complex
$$0 \r \bar{K}^{\mathrm{M}}_n(A) \r
\bigoplus_{x \in A^{(0)}} \KM_i\big(k(x)\big)
\stackrel{\partial}{\longrightarrow} \bigoplus_{x \in A^{(1)}}
\KM_{i-1}\big(k(x)\big) \stackrel{\partial}{\longrightarrow} \ldots,$$
where $A^{(r)}$ is the set of codimension-$r$ points in $\mathrm{Spec}
\ A$ and $\kappa(x)$ is the residue field at $x$. In order to
determine all operations $\bar{K}^{\mathrm{M}}_i/p \r
\bar{K}^{\mathrm{M}}_*/p$ over a ``nice'' ring $A$, we will be
concerned with the exactness of that complex. Let $A$ be an
essentially smooth semi-local $k$-algebra, where by essentially smooth
we mean that $A$ is a localization of a smooth affine $k$-algebra.
Gabber (unpublished), Elbaz-Vincent and Mueller-Stach have established
the exactness of the complex $$ \bar{K}^{\mathrm{M}}_n(A) \r
\bigoplus_{x \in A^{(0)}} \KM_i\big(k(x)\big)
\stackrel{\partial}{\longrightarrow} \bigoplus_{x \in A^{(1)}}
\KM_{i-1}\big(k(x)\big) \stackrel{\partial}{\longrightarrow} \ldots$$
in the case when $A$ has infinite residue field (see \cite{EVMS} and
\cite{KMS}). This last condition has been removed by Kerz in
\cite{Kerz2}. Moreover, in \cite{Kerz}, Kerz shows that the Gersten
complex is also exact at the first place when $A$ has infinite residue
field. All in all,

\begin{theorem}[Gabber, Elbaz-Vincent, Mueller-Stach, Kerz]
  \label{KMS} 
  Let $A$ be an essentially smooth semi-local algebra over a field $k$
  with quotient field $F$. If $A$ has infinite residue field, the
  Gersten complex is exact, and in particular
$$\bar{K}^{\mathrm{M}}_n(A) = \im \Big( A^\times \otimes \ldots \otimes
A^\times \longrightarrow \KM_n(F) \Big)= \ker \Big(\KM_n(F)
\stackrel{\partial}{\longrightarrow} \bigoplus_{x \in A^{(1)}}
\KM_{n-1}\big(k(x)\big) \Big).$$ Without assuming $A$ has infinite
residue field, $$\im \Big( A^\times \otimes \ldots \otimes A^\times
\longrightarrow \KM_n(F) \Big)= \ker \Big(\KM_n(F)
\stackrel{\partial}{\longrightarrow} \bigoplus_{x \in A^{(1)}}
\KM_{n-1}\big(k(x)\big) \Big).$$
\end{theorem}

From now on, $k$ is an infinite field and $A$ is a fixed essentially
smooth semi-local $k$-algebra. Let's denote by $\mathsf{C}_A$ the
category of rings over $A$ (i.e. the rings $R$ with a morphism $R\r
A$) with morphisms compatible with the structure maps to $A$. Note
that fields containing $A$ are objects in the category $\mathsf{C}_A$.
The map $\bar{K}^{\mathrm{M}}_n/p : \mathsf{C}_A \r \mathsf{Sets}$ is
a functor, and an operation over the regular semi-local domain $A$ is
a natural transformation from the functor $\bar{K}^{\mathrm{M}}_n/p :
\mathsf{C}_A \r \mathsf{Sets}$ to the functor
$\bar{K}^{\mathrm{M}}_*/p : \mathsf{C}_A \r \mathbf{F}_p -
\mathsf{Algebras}$. Under these assumptions, all the results
concerning fields translate to the case of essentially smooth
semi-local $k$-algebras and it is possible to describe all such
operations. \medskip

\begin{theorem}\label{oplocal} Let $k$ be an infinite field and $A$ be
  an essentially smooth semi-local $k$-algebra.  Operations
  $\bar{K}^{\mathrm{M}}_i/p \r \bar{K}^{\mathrm{M}}_*/p$ over $A$
  are in the $\KM_*(A)/p$-span of divided powers.
\end{theorem}

\noindent \textbf{Proof.}  Indeed, by Theorem \ref{KMS}, the functors
$ \bar{K}^{\mathrm{M}}_*/p$ and $\KM_*/p$ agree on essentially smooth
semi-local $k$-algebras. The proof of Theorem \ref{oprings} applied to
$\mathrm{Spec} \ A$ shows that an operation over $A$ must be in the
span of divided power operations as in Proposition \ref{dprings}. It
remains to say that such operations do exist. For this purpose, it is
enough to check that $\{b,b\}=\{-1,b\} \in \bar{K}^{\mathrm{M}}_*(B)$
for all $k$-algebra $B$ and all $b \in B^\times$. This is the case
because $k$ is infinite.  Indeed, the relation $\{b,-b\}=0 \in
\bar{K}^{\mathrm{M}}_*(B)$ holds for all $b \in B^\times$ whenever the
field $k$ is infinite, see e.g.  \cite{NS}.\qed \medskip

\section{Operations from Milnor $K$-theory to a cycle module for
  fields}

The notion of cycle module has been defined and thoroughly studied by
Rost in \cite{Rost}. A cycle module is a $\Z$-graded functor $M_* :
\mathsf{Fields} \r \mathsf{Abgroups}$ equipped with four structural
data and satisfying certain rules and axioms. In particular, for all
field $k$, $M_*(k)$ is a left $\KM_*(k)$-module such that the product
respects the grading. It is also a right $\KM_*(k)$-module in the
following way : if $\rho \in \KM_n(k)$ and $x \in M_m(k)$, then $x
\cdot \rho = (-1)^{mn} \rho \cdot x$. Examples of cycle modules are
given by Galois cohomology, Milnor $K$-theory and Quillen $K'$-theory,
see Theorem 1.4 and Remark 2.5 of \textit{loc. cit.}. These examples
are actually cycle modules with a ring structure \cite[Def 1.2]{Rost}.

Our results in this section applied to Quillen's $K$-theory of fields
say that operations mod $p$ from Milnor $K$-theory to Quillen
$K$-theory are spanned by divided power operations. Precisely, let
$k_0$ be a field and $p$ be a prime number. We are interested in
operations $\KM_i/p \r K_*/p$ over the field $k_0$, that is, natural
transformations from the Milnor $K$-theory functor $\KM_i/p :
\mathsf{Fields}_{/k_0} \r \mathsf{Sets}$ to the Quillen $K$-theory
functor $K_*/p : \mathsf{Fields}_{/k_0} \r
\mathbf{F}_p-\mathsf{Algebras}$. There is a natural transformation
$\KM_* \r K_*$ for fields induced by cup-product, which is identity in
degrees $0$ and $1$. A divided power $\gamma_n : \KM_i(k) \r
K_{ni}(k)$ is the composition of the divided power $\gamma_n :
\KM_i(k) \r \KM_{ni}(k)$ on Milnor $K$-theory with the natural map
$\KM_{ni}(k) \r K_{ni}(k)$.  Divided powers are well-defined in the
same cases as divided powers in Milnor K-theory, e.g. mod $p$ for $p$
odd and $i$ even. Note that composing with the natural map
$\KM_{ni}(k) \r K_{ni}(k)$ amounts to multiply by $1$ when $K_*(k)$ is
seen as a $\KM_*(k)$-module. Obviously, divided powers, when
well-defined, are indeed operations. Theorem \ref{cyclemod} say that
divided power operations span as a $K_*(k_0)$-module the algebra of
operations $\KM_i/p \r K_*/p$ over $k_0$.

Some properties of cycle modules of most interest to us are the following

\begin{proposition}[Homotopy property for $\mathbf{A}^1$]
  \label{homotopy} Let $k$ be a field. There is a short exact sequence
  $$0 \longrightarrow M_*(k) \stackrel{\iota}{\longrightarrow}
  M_*(k(t)) \stackrel{\partial}{\longrightarrow} \bigoplus_{P \in
    (\mathbf{A}^1_k)^{(1)}} M_{*-1}(\kappa(P)) \longrightarrow 0.$$
  The map $\iota$ is induced by the inclusion of fields $k
  \hookrightarrow k(t)$ and $\partial$ is the sum of the residue maps
  at the closed points of the affine line $\mathbf{A}^1_k$.
\end{proposition}

\noindent \textbf{Proof.} This is proposition 2.2. of \cite{Rost}. \qed

\begin{corollary}
  If $\partial_0 : M_*(k(t)) \r M_*(k)$ is the residue map
  at $0$, then for any $x \in M_*(k)$, we have the formula
  $$\partial_0(\{t\}\cdot x)=x.$$ In particular, the map $M_{*-1}(k)
  \r M_{*}(k(t)), \ x \mapsto x\cdot \{t\}$ is injective.
\end{corollary}

\noindent \textbf{Proof.} By the previous proposition, we have
$\partial_0(x)=0$. The rule R3f of \cite{Rost} then gives
$\partial_0(\{t\}\cdot x) = \partial_0(\{t\}) \cdot s_t(x) =
s_t(x)=x$.  \qed

\subsection{Operations mod $p$}

Let $k_0$ be a field and $p$ be a prime number. We are interested in
operations $\KM_i/p \r M_*/p$ over the field $k_0$, that is, natural
transformations from the Milnor $K$-theory functor $\KM_i/p :
\mathsf{Fields}_{/k_0} \r \mathsf{Sets}$ to the cycle module functor
$M_*/p : \mathsf{Fields}_{/k_0} \r \mathbf{F}_p-\mathsf{Algebras}$. A
divided power is a map $a\cdot \gamma_n : \KM_i(k) \r M_*(k)$, where
$a \in M_*(k)$ and $\gamma_n$ is the divided power defined on Milnor
$K$-theory. Theorems \ref{main} and \ref{main2} generalize to cycle
modules with ring structure.

\begin{theorem} \label{cyclemod} Let $k_0$ be any field and $p$ be a
  prime number.  Suppose $M_*$ is a cycle module with ring structure.
  The algebra of operations $\KM_i/p \r M_*/p$ over $k_0$ is
\begin{itemize}
\item If $i=0$, the free $M_*(k_0)/p$-module of rank $p$ of
  functions $\mathbf{F}_p \r M_*(k_0)/p$.
\item If $i=1$, the free $M_*(k_0)/p$-module of rank $2$ generated by
  $\gamma_0$ and $\gamma_1$.
\item If $i\geq 1$ odd and $p$ odd, the free $M_*(k_0)/p$-module of
  rank $2$ generated by $\gamma_0$ and $\gamma_1$.
\item If $i \geq 2$ even and $p$ odd, the free $M_*(k_0)/p$-module
  $$\bigoplus_{n\geq 0} M_*(k_0)/p \cdot \gamma_n.$$
\item If $i \geq 2$ and $p=2$, the $M_*(k_0)/2$-module $$M_*(k_0)/2
  \cdot \gamma_0 \oplus M_*(k_0)/2 \cdot \gamma_1 \oplus
  \bigoplus_{n\geq 2} \mathrm{Ker}(\tau_i) \cdot \gamma_n,$$ where
  $\tau_i$ is the map $M_*(k_0)/2 \r M_*(k_0)/2, \ x \mapsto
  \{-1\}^{i-1} \cdot x$.
\end{itemize}
\end{theorem}

\noindent \textbf{Proof.}  The main ingredients used in the proof of
this theorem in the case of operations from Milnor $K$-theory to
itself were the following. We first determined the operations $\KM_1/p
\r \KM_*/p$ over $k_0$, and showed that they were of the form $x
\mapsto a\cdot x +b$.  From there, it was easy to determine the
operations $\big(\KM_1/p\big)^r \r \KM_*/p$. In order to determine the
operations $\KM_i/p \r \KM_*/p$ for $i\geq 2$, we used the fact that
Milnor $K$-theory of a field is generated as a ring by the degree $1$
elements and also the important fact that for $t$ transcendental over
$k_0$ the map $\KM_{n-1}(k_0) \r \KM_{n}(k_0(t)),\ x \mapsto x\cdot
\{t\}$ is injective.

The case $i=0$ is easy.

Let's first deal with the case $i=1$. Let $\varphi : \KM_1 \r M_*$ be
an operation over $k_0$. Then, $\varphi$ is determined by the image of
a transcendental element $t$ over $k_0$. This is because the map
$M_*(k) \r M_*(k(t))$ induced by the transcendental extension $k(t)/k$
is injective by proposition \ref{homotopy}. We claim that the element
$\varphi(\{t\}) \in M_*(k_0(t))/p$ has residue $0$ for all residue
maps corresponding to closed points in $\A_{k_0}^1 \backslash \{0\}$.
The proof is essentially the same as step 2 of the proof of Theorem
\ref{theorem} using the homotopy property of proposition
\ref{homotopy} in place of Milnor's exact sequence.

So, $\varphi(\{t\})$ is unramified outside $0$.  Let $b =
\partial_0\big(\varphi(\{t\})\big)$, then thanks to proposition
\ref{homotopy} and its corollary, we get the existence of $a \in
M_{*}(k_0)/p$ such that $$\varphi(\{t\}) = a + \{t\} \cdot b \in
M_*(k_0(t))/p.$$ This defines an operation as one can easily check.

From there, it is routine to check (cf. Step 3 of the proof of theorem
\ref{theorem}) that given an operation $\varphi : \big(\KM_1/p\big)^r
\r M_*/p$, there exist elements $\lambda_{i_1,\ldots,i_s} \in
M_*(k_0)/p$ such that for all field extension $k/k_0$ and all
$r$-tuple $(\{a_1\},\ldots,\{a_r\}) \in \big(\KM_1(k)/p \big)^{r}$,
$$\varphi(\{a_1\},\ldots,\{a_r\}) = \sum_{1\leq i_1 < \ldots < i_s \leq r}
\lambda_{i_1,\ldots,i_s} \cdot \{a_{i_1},\ldots,a_{i_s}\}.$$

It remains to prove the cases for which $i \geq 2$. In each cases of
the theorem, divided powers are well-defined and do define operations.
It is thus enough to prove that an operation $\varphi : \KM_i/p \r
M_*/p$ over $k_0$ must be of the form stated. A close look at the
proof of proposition \ref{main1}, as well as the proofs of
propositions \ref{odd} and \ref{opmod2}, shows that the property asked
to the functor $M_*$ for the proofs to translate \textit{mutatis
  mutandis} to the case of operations $\KM_i/p \r M_*/p$ is that the
map $M_{*-1}(k) \r M_{*}(k(t)), \ x \mapsto x\cdot \{t\}$ is
injective. This is the object of the corollary to proposition
\ref{homotopy}. \qed

\subsection{Integral operations}

Let $k$ be a field with discrete valuation $v$, and let $\pi$ be a
uniformizer for $v$. Given a cycle module $M_*$, there is a
specialization map $s_\pi : M_*(k) \r M_*(\kappa(v))$ defined by
$$s_\pi(x) = \partial_v(\{-\pi\}\cdot x).$$

\begin{definition}
  Let $k_0$ be any field and $K$ an extension of $k_0$ endowed with a
  discrete valuation $v$ such that its valuation ring $R = \{x \in K,
  v(x) \geq 0\}$ contains $k_0$, so that the residue field $\kappa$ is
  an extension of $k_0$. We say that specialization maps commute with
  an operation $\varphi : \KM_i \longrightarrow M_*$ over $k_0$ if for
  any extension $K/k_0$ as above, we have a commutative diagram
  \begin{center} $ \xymatrix{\KM_i(K)/p \ar[r]^{
        \ \varphi} \ar[d]_{s_\pi} & M_*(K)/p \ar[d]_{s_\pi} \\
      \KM_i(\kappa)/p \ar[r]^{ \ \varphi} & M_*(\kappa)/p }$
  \end{center} where $\pi$ is any uniformizer for the valuation $v$.
\end{definition}

\begin{theorem} Let $k_0$ be a field and suppose $M_*$ is a cycle
  module with ring structure. The algebra of operations $\KM_i \r M_*$
  over $k_0$ commuting with specialization maps is

\begin{itemize}
\item If $i=0$, the $M_*(k_0)$-module of functions $\Z \r M_*(k_0)$.
\item If $i=1$, the free $M_*(k_0)$-module generated by $\gamma_0$ and
  $\gamma_1$.
\item If $i$ is even $\geq 2$, the $M_*(k_0)$-module $$M_*(k_0) \oplus
  M_*(k_0)\cdot\mathrm{id} \oplus \bigoplus_{n\geq 2}
  \mathrm{Ker}(\tau_i)\cdot \gamma_n.$$
\item If $i$ is odd $\geq 2$, the $M_*(k_0)$-module $$M_*(k_0) \oplus
  M_*(k_0)\cdot\mathrm{id} \oplus \bigoplus_{n\geq 2}
  {}_2\mathrm{Ker}(\tau_i)\cdot \gamma_n.$$
\end{itemize} Here, $\tau_i$ is the map $M_*(k_0) \r M_*(k_0), \ x
\mapsto \{-1\}^{i-1} \cdot x$.
\end{theorem}

\noindent \textbf{Proof.} Same as for Milnor $K$-theory. See section 3.4.\qed

\begin{theorem}
  Let $k_0$ be a field and suppose $M_*$ is a cycle module with ring
  structure. The algebra of additive operations $\varphi : \KM_i \r
  M_*$ over $k_0$ is the free $M_*(k_0)$-module generated by the
  identity.
\end{theorem}

\noindent \textbf{Proof.} Same as for Milnor $K$-theory. See section 3.5. \qed

\begin{corollary}
  A natural map from $\KM_i : \mathsf{Fields} \r \mathsf{Sets}$ to
  $K_i : \mathsf{Fields} \r \mathsf{AbGroups}$ is an integral multiple
  of the usual natural map.
\end{corollary}

\subsection{Operations for smooth schemes}

As in the case of Milnor $K$-theory, the operations $\KM_i/p \r
H^0(-,M_*)$ for smooth schemes over $k$ are spanned as a
$M_*(k)/p$-module by divided powers.

\begin{theorem} Let $k$ be any field, $p$ be a prime number and $X$ be
  a smooth scheme over $k$. Operations $\KM_i/p \r H^0(-,M_*/p)$ over
  the smooth $k$-scheme $X$ are spanned as a $H^0(X,M_*/p)$-module by
  the divided powers.
\end{theorem}

\noindent \textbf{Proof.} Same as for Milnor $K$-theory. See section 4.1. \qed

\begin{small}

\def\cprime{$'$}
  \medskip

\textsc{Department of Pure Mathematics and Mathematical Statistics,
  University of Cambridge, Wilberforce Road, Cambridge, CB3 0WB,
  United Kingdom}
 \end{small}

\textit{e-mail :}  \texttt{C.Vial@dpmms.cam.ac.uk}


\begin{thebibliography}{10}

\bibitem{BassTate}
H.~Bass and J.~Tate.
\newblock The {M}ilnor ring of a global field.
\newblock In {\em Algebraic $K$-theory, II: ``Classical'' algebraic $K$-theory
  and connections with arithmetic (Proc. Conf., Seattle, Wash., Battelle
  Memorial Inst., 1972)}, pages 349--446. Lecture Notes in Math., Vol. 342.
  Springer, Berlin, 1973.

\bibitem{Becher}
K.~J. Becher.
\newblock Virtuelle {F}ormen.
\newblock In {\em Mathematisches Institut, Georg-August-Universit\"at
  G\"ottingen: Seminars 2003/2004}, pages 143--150. Universit\"atsdrucke
  G\"ottingen, G\"ottingen, 2004.

\bibitem{Becher2}
K.~J. Becher and D.~W. Hoffmann.
\newblock Symbol lengths in {M}ilnor {$K$}-theory.
\newblock {\em Homology Homotopy Appl.}, 6(1):17--31 (electronic), 2004.

\bibitem{BO}
Pierre Berthelot and Arthur Ogus.
\newblock {\em Notes on crystalline cohomology}.
\newblock Princeton University Press, Princeton, N.J., 1978.

\bibitem{CT}
J.-L. Colliot-Th{\'e}l{\`e}ne.
\newblock Birational invariants, purity and the {G}ersten conjecture.
\newblock In {\em $K$-theory and algebraic geometry: connections with quadratic
  forms and division algebras (Santa Barbara, CA, 1992)}, volume~58 of {\em
  Proc. Sympos. Pure Math.}, pages 1--64. Amer. Math. Soc., Providence, RI,
  1995.

\bibitem{EVMS}
Philippe Elbaz-Vincent and Stefan M{\"u}ller-Stach.
\newblock Milnor {$K$}-theory of rings, higher {C}how groups and applications.
\newblock {\em Invent. Math.}, 148(1):177--206, 2002.

\bibitem{EL1}
Richard Elman and T.~Y. Lam.
\newblock Quadratic forms over formally real fields and pythagorean fields.
\newblock {\em Amer. J. Math.}, 94:1155--1194, 1972.

\bibitem{EL2}
Richard Elman and T.~Y. Lam.
\newblock Classification theorems for quadratic forms over fields.
\newblock {\em Comment. Math. Helv.}, 49:373--381, 1974.

\bibitem{Epstein}
D.~B.~A. Epstein.
\newblock Steenrod operations in homological algebra.
\newblock {\em Invent. Math.}, 1:152--208, 1966.

\bibitem{Galois}
Skip Garibaldi, Alexander Merkurjev, and Jean-Pierre Serre.
\newblock {\em Cohomological invariants in {G}alois cohomology}, volume~28 of
  {\em University Lecture Series}.
\newblock American Mathematical Society, Providence, RI, 2003.

\bibitem{Gille}
Philippe Gille and Tam{\'a}s Szamuely.
\newblock {\em Central simple algebras and {G}alois cohomology}, volume 101 of
  {\em Cambridge Studies in Advanced Mathematics}.
\newblock Cambridge University Press, Cambridge, 2006.

\bibitem{Izhboldin}
O.~Izhboldin.
\newblock On {$p$}-torsion in {$K\sp M\sb *$} for fields of characteristic
  {$p$}.
\newblock In {\em Algebraic $K$-theory}, volume~4 of {\em Adv. Soviet Math.},
  pages 129--144. Amer. Math. Soc., Providence, RI, 1991.

\bibitem{KahnSW}
Bruno Kahn.
\newblock Classes de {S}tiefel-{W}hitney de formes quadratiques et de
  repr\'esentations galoisiennes r\'eelles.
\newblock {\em Invent. Math.}, 78(2):223--256, 1984.

\bibitem{Kafields}
Bruno Kahn.
\newblock Comparison of some field invariants.
\newblock {\em J. Algebra}, 232(2):485--492, 2000.

\bibitem{Kerz}
Moritz Kerz.
\newblock {Milnor $K$-theory of local rings with finite residue fields }.
\newblock {Preprint, September 11, 2007, K-theory Preprint Archives,
  http://www.math.uiuc.edu/K-theory/0865/}.

\bibitem{Kerz2}
Moritz Kerz.
\newblock The {G}ersten conjecture for {M}ilnor {$K$}-theory.
\newblock {\em Invent. Math.}, 175(1):1--33, 2009.

\bibitem{KMS}
Moritz Kerz and Stefan M{\"u}ller-Stach.
\newblock The {M}ilnor-{C}how homomorphism revisited.
\newblock {\em $K$-Theory}, 38(1):49--58, 2007.

\bibitem{Lam}
T.~Y. Lam.
\newblock {\em Introduction to quadratic forms over fields}, volume~67 of {\em
  Graduate Studies in Mathematics}.
\newblock American Mathematical Society, Providence, RI, 2005.

\bibitem{MacDonald}
Mark MacDonald.
\newblock Cohomological invariants of odd degree {J}ordan algebras.
\newblock {\em Math. Proc. Cam. Phil. Soc., to appear.}

\bibitem{MS}
A.~S. Merkurjev and A.~A. Suslin.
\newblock {$K$}-cohomology of {S}everi-{B}rauer varieties and the norm residue
  homomorphism.
\newblock {\em Izv. Akad. Nauk SSSR Ser. Mat.}, 46(5):1011--1046, 1135--1136,
  1982.

\bibitem{Milnor}
John Milnor.
\newblock Algebraic {$K$}-theory and quadratic forms.
\newblock {\em Invent. Math.}, 9:318--344, 1969/1970.

\bibitem{NS}
Yu.~P. Nesterenko and A.~A. Suslin.
\newblock Homology of the general linear group over a local ring, and
  {M}ilnor's {$K$}-theory.
\newblock {\em Izv. Akad. Nauk SSSR Ser. Mat.}, 53(1):121--146, 1989.

\bibitem{Revoy}
Philippe Revoy.
\newblock Formes altern\'ees et puissances divis\'ees.
\newblock In {\em S\'eminaire {P}. {D}ubreil, 26e ann\'ee (1972/73),
  {A}lg\`ebre, {E}xp. {N}o. 8}, page~10. Secr\'etariat Math\'ematique, Paris,
  1973.

\bibitem{Rost}
Markus Rost.
\newblock Chow groups with coefficients.
\newblock {\em Doc. Math.}, 1:No. 16, 319--393 (electronic), 1996.

\bibitem{Totaro}
Burt Totaro.
\newblock Milnor {$K$}-theory is the simplest part of algebraic {$K$}-theory.
\newblock {\em $K$-Theory}, 6(2):177--189, 1992.

\bibitem{Voevodsky2}
Vladimir Voevodsky.
\newblock {On motivic cohomology with $\mathbf{Z}/l$-coefficients}.
\newblock {Preprint, June 16, 2003, K-theory Preprint Archives,
  http://www.math.uiuc.edu/K-theory/0639/}.

\bibitem{Voevodsky}
Vladimir Voevodsky.
\newblock Motivic cohomology with {${\bf Z}/2$}-coefficients.
\newblock {\em Publ. Math. Inst. Hautes \'Etudes Sci.}, (98):59--104, 2003.

\bibitem{Voevodsky1}
Vladimir Voevodsky.
\newblock Reduced power operations in motivic cohomology.
\newblock {\em Publ. Math. Inst. Hautes \'Etudes Sci.}, (98):1--57, 2003.

\bibitem{Weibel}
Charles~A. Weibel.
\newblock {Patching the Norm Residue Isomorphism Theorem }.
\newblock {Preprint, May 23, 2007, K-theory Preprint Archives,
  http://www.math.uiuc.edu/K-theory/0844/}.

\end{thebibliography}
\end{document}